\documentclass[]{article}

\usepackage{amsmath,amssymb,amsfonts,xfrac,MnSymbol}
\usepackage{amsthm}
\usepackage{cite}
\usepackage{hyperref}
\usepackage{todonotes}
\usepackage{enumitem}
\usepackage{graphicx}
\usepackage{color}
\usepackage{import}

\numberwithin{equation}{section}

\newtheorem{theorem}{Theorem}[section]

\newtheorem{proposition}[theorem]{Proposition}
\newtheorem{lemma}[theorem]{Lemma}
\newtheorem{corollary}[theorem]{Corollary}

\newtheorem*{claim*}{Claim}

\newtheorem{theoremA}{Theorem}

\theoremstyle{definition}
\newtheorem{definition}[theorem]{Definition}

\newtheorem{remark}[theorem]{Remark}

\newtheorem*{notation*}{Notation}

\newcommand{\bbN}{\mathbb{N}}

\newcommand{\actson}{\curvearrowright}

\newcommand{\epi}{\twoheadrightarrow}
\newcommand{\into}{\hookrightarrow}

\newcommand{\normalin}{\lhd}

\newcommand{\ii}{^{-1}}

\newcommand{\gen}[1]{\left< #1 \right>}
\newcommand{\ngen}[1]{\left\llangle #1 \right\rrangle}

\newcommand{\qc}{quasiconvex}

\newcommand{\qg}[1]{(#1)-quasigeodesic}

\newcommand{\torsionscc}{C_1(\epsilon,\mu,\lambda,c,\rho)}

\newcommand{\torsionsccK}[1]{C_1(\epsilon,\mu,\lambda,c,\rho,#1)}

\newcommand{\fresh}[1]{#1^{\natural}}

\DeclareMathOperator{\Aut}{Aut}

\DeclareMathOperator{\Core}{Core}

\DeclareMathOperator{\diam}{diam}

\title{Invariable generation does not pass to finite index subgroups}
\author{Gil Goffer and Nir Lazarovich}
\date{}

\begin{document}
\maketitle

\begin{abstract}
    Using small cancellation methods, we show that the property \emph{invariable generation} does not pass to finite index subgroups, answering questions of Wiegold \cite{wiegold1977transitive} and Kantor-Lubotzky-Shalev \cite{kantor2011invariable}.
    We further show that a finitely generated group that is invariably generated is not necessarily finitely invariably generated, answering a question of Cox \cite{cox2020invariable}. The same results were also obtained independently by Minasyan \cite{minasyan2020invariable}.
\end{abstract}

\section{Introduction}

\begin{definition}[Dixon \cite{dixon1992random}] Let $G$ be a group. A subset $S\subseteq G$ \emph{invariably generates} $G$ if for every function $S\to  G, s\mapsto g_s$, the set of conjugates $\{s^{g_s}|s\in S\}$ generates $G$.

A group $G$ is \emph{invariably generated} (or \emph{IG}) if it has an invariably generating set. That is, if $G$ invariably generates itself. A group $G$ is \emph{finitely invariably generated} (or \emph{FIG}), if it has a finite invariably generating set.
\end{definition}

%The notion was first introduced by Dixon \cite{dixon1992random}, for finite groups, \cite{dixon1992random}, motivated by Chebotarev’s Density Theorem \cite{tschebotareff1926bestimmung}.
Dixon's original definition referred to finite groups. %, motivated by Chebotarev’s Density Theorem \cite{tschebotareff1926bestimmung}.
However, an equivalent definition was previously studied by Wiegold in the context of general (finite or infinite) groups  \cite{wiegold1976transitive}.
Kantor, Lubotzky and Shalev \cite{kantor2011invariable} were the first to consider Dixon's definition for infinite groups, and to notice that it coincides with Wiegold's definition.

It is shown in \cite{wiegold1976transitive, kantor2011invariable} %\todo{N: is this the right reference?}
that the classes of IG groups and FIG groups are closed under extensions and include all finite groups. It follows that a group with a finite index normal IG (resp. FIG) subgroup is IG (resp. FIG). The following slight generalization is probably known to experts, yet we include a proof of this theorem in Section \ref{sec: the other direction}.

\begin{theoremA}\label{thm: other direction (intro)}
    A group containing a finite index IG (resp. FIG) subgroup is IG (resp. FIG).
\end{theoremA}

In contrast, we prove the following theorem, answering questions of Wiegold \cite{wiegold1977transitive} and Kantor-Lubotzky-Shalev \cite{kantor2011invariable}.

\begin{theoremA} \label{thm: main}
	There exists a FIG group with an index 2 non-IG subgroup.
\end{theoremA}

In the context of topological groups, it was shown in \cite{kantor2011invariable} that a topologically finitely generated group that is topologically invariably generated is not necessarily finitely invariably generated.
We therefore find it relevant to state the following theorem, answering a question of Cox \cite{cox2020invariable}.

\begin{theoremA}\label{thm: there exists fg IG not FIG group}
There exists a finitely generated group that is invariably generated, but not finitely invariably generated.
\end{theoremA}

The proofs of Theorem \ref{thm: main} and Theorem \ref{thm: there exists fg IG not FIG group} rely on an iterative small cancellation construction. The same results were obtained independently by \mbox{Minasyan} \cite{minasyan2020invariable} using similar methods.

%We remark that from the historical point of view, the notion of invariable generation was initially considered for finite groups \cite{dixon1992random,kantor2011invariablefinite}, motivated by Chebotarev’s Density Theorem \cite{tschebotareff1926bestimmung}. Further works in this context consider invariable generation of the symmetric groups $S_n$ and other classical finite groups \cite{eberhard_kevin_ford2017invariable,pemantle_peres_rivin2016four,mckemmie2019invariable}.
%Kantor, Lubotzky and Shalev were the first ones to extend Dixon's definition to to infinite groups \cite{kantor2011invariable}, and to notice that it coincides with the (independently developed) Wiegold's definition \cite{wiegold1976transitive}. %(Items \ref{equivalent definitions - conjugacy complete subgroup} and \ref{equivalent definitions - transitive action} in Lemma \ref{lem: equivalent definitions of IG}).
Invariable generation was studied for various groups and classes of groups, including
symmetric groups  \cite{dixon1992random,eberhard_kevin_ford2017invariable,pemantle_peres_rivin2016four},
finite groups \cite{detomi_lucchini2015invariable,mckemmie2019invariable,kantor2011invariablefinite},
wreath products \cite{lucchini2017invariable,cox2020invariable},
the Thomspon groups \cite{gelander_golan_Juschenko2017invariable},
convergence groups \cite{gelander2015convergence},
linear groups \cite{kantor2011invariable,gelander_meiri2017congruence} and topological groups \cite{goffer2018few,kantor2011invariable}.

%of finite groups was studied in \cite{kantor2011invariablefinite,eberhard_kevin_ford2017invariable,pemantle_peres_rivin2016four,mckemmie2019invariable},
%and of infinite (discrete and topological) groups in   \cite{gelander2015convergence,gelander_meiri2017congruence,gelander_golan_Juschenko2017invariable,goffer2018few}.
%Finitely invariable generation for infinite groups was recently studies in \cite{detomi_lucchini2015invariable,lucchini2017invariable,cox2020invariable}.
%Under this scope, we find Theorems \ref{thm: main} and \ref{thm: there exists fg IG not FIG group} very relevant*\todo{G: replace the word 'relevant' by something. Maybe 'applicable'?}.

\emph{Organization of the paper:}
In Section \ref{sec: the other direction} we include the proof of Theorem \ref{thm: other direction (intro)}.
In Section \ref{sec: toolbox} we give a brief statement of the tools used in the proofs of Theorems \ref{thm: main} and \ref{thm: there exists fg IG not FIG group}. In Section \ref{sec: constructing G} we prove Theorem \ref{thm: main}. In Section \ref{sec: constructing fg IG not FIG group} we prove Theorem \ref{thm: there exists fg IG not FIG group}. In Section \ref{sec: preliminaries on hyp groups and sc} we give the main definitions for small cancellation theory of hyperbolic groups following Ol'shanskii \cite{olshanskii1993residualing}. In Section \ref{sec: finding sc words}, we show that one can find small cancellation words with specific properties, and prove the main lemmas of Section \ref{sec: toolbox}. Section \ref{sec: hexagon} is devoted to the hexagon property which is an ingredient of the proof of Theorem \ref{thm: main}.

\section{Proof of Theorem \ref{thm: other direction (intro)}}
\label{sec: the other direction}

\begin{definition}
	Let $G$ be a group, and $S\subseteq G$ a subset. A subgroup $H \le G$ is \emph{$S$-conjugacy complete} if it intersects the conjugacy classes of all elements of $S$.
	
	When $S=G$ we say that $H$ is \emph{conjugacy complete}.
\end{definition}

It was observed in \cite{kantor2011invariable} that the following are equivalent definitions of IG.
\begin{lemma}\label{lem: equivalent definitions of IG}
	Let $G$ be a group, and $S\subseteq G$ a subset.
	The following are equivalent:
	\begin{enumerate}
	    \item $S$ invariably generates $G$ \label{equivalent definitions - IG}
	    \item $G$ does not contain a proper $S$-conjugacy complete subgroup.\label{equivalent definitions - conjugacy complete subgroup}
	    \item Every non-trivial transitive action $G\actson X$ has an element $s\in S$ without fixed points.\label{equivalent definitions - transitive action}
	\end{enumerate}
\end{lemma}

Wiegold \cite{wiegold1976transitive}
proved that the class of IG groups is closed under extensions, in fact the following slightly stronger result holds.

\begin{proposition}\label{prop: other direction}
Let $G$ be a group,
%and let $N\normalin G$, and $N\le H \le G$ be subgroups.
$N\le H \le G$ be subgroups and $N\normalin G$.
Let $S\subseteq H$ and $S'\subseteq G$. If $H$ is invariably generated by $S$ and $G/N$ is invariably generated by $S'$ then $G$ is invariably generated by $S\cup S'$.
\end{proposition}

\begin{proof}
    Let $G\actson X$ be a transitive action on a set with $|X|\ge 2$. We want to find an element of $S\cup S'$ which acts without fixed points on $X$.

    Since $N$ is normal, we know that $G/N\actson X/N$. If $|X/N|\ge 2$ then since $G/N$ is invariably generated by $S'$, there exists an element $s' \in S'$ that acts without fixed points on $X/N$ and hence also on $X$. If $|X/N|=1$, then $N$, and hence $H$, act transitively on $X$ and since $H$ is invariably generated by $S$, there is an element $s\in S$ which acts without fixed points on $X$.
\end{proof}

In particular, we can deduce Theorem \ref{thm: other direction (intro)}.
\begin{proof}[Proof of Theorem \ref{thm: other direction (intro)}]
If $H$ is a finite index IG (resp. FIG) subgroup of $G$, then $N=\Core_G(H)=\bigcap_{g\in G} g\ii Hg$ is of finite index in $G$. Since every finite group is FIG, we get that $N\le H\le G$ satisfy the assumptions of Proposition \ref{prop: other direction} which implies that $G$ is IG (resp. FIG).
\end{proof}

\section{Toolbox}\label{sec: toolbox}

In this section we describe the toolbox for the main constructions.
Since the main constructions are based on small cancellation quotients and HNN extensions, we summarize in this section the main relevant lemmas regarding these two topics.
We believe that a reader who is familiar with small cancellation theory would feel fairly comfortable with these lemmas, whose proofs follow standard techniques. We therefore postpone their proofs to later sections.

Throughout the rest of the paper we assume familiarity with notions in hyperbolic group theory (cf. for example \cite{bridson2013metric,olshanskii1993residualing,gromov1987hyperbolic,ghys1990espaces}).

\subsection{Small cancellation quotients}

We use the small cancellation theory developed by Ol'shanskii \cite{olshanskii1993residualing} for hyperbolic groups, which we outline more precisely in Section \ref{sec: preliminaries on hyp groups and sc}.

Roughly speaking, we say that a set of quasigeodesic words $\mathcal{R}$ in a hyperbolic group satisfies small cancellation if whenever two words in $\mathcal{R}$ fellow-travel, they do so for a small proportion of their lengths.
Similarly, we say that a set of quasigeodesic words $\mathcal{R}$ has small overlap with another set of quasigeodesic words $\mathcal{K}$, if whenever a word in $\mathcal{R}$ fellow-travels a word in $\mathcal{K}$, it does so for a small proportion of its length.

\begin{lemma} \label{lem: existence of small enough sc words}
    Let $G$ be a torsion-free hyperbolic group, and let $H,K_1,\ldots,K_n$ be quasiconvex subgroups of $G$. If $H$ is non-elementary and non-commensurable\footnote{we use the term ``commensurable'' to refer to the equivalence of subgroups up to conjugation and passing to finite index. That is, two subgroups $H,H'$ in $G$ are \emph{commensurable} if there exists $g\in G$ such that $H^g\cap H'$ has finite index in both $H^g$ and $H'$. Similarly, $H$ is \emph{commensurable into} $H'$ if there exists $g\in G$ such that $H^g \cap H'$ has finite index in $H^g$.} %\todo{N:what do you think of this footnote?}
    into $K_1,\ldots,K_n$.
    Then for every $m$ there exists a subset of $m$ words $\mathcal{R}=\{w_1,\ldots,w_m\}\subseteq H$ with arbitrarily small cancellation and arbitrarily small overlap with $K_1,\ldots,K_n$.%\todo{N: add E(w)=w here? or should we leave it as a remark after the proof of this lemma?}

    If moreover $G$ has an involution $\phi$ which exchanges two non-commensurable\footnote{elements are \emph{commensurable} if they generate cyclic subgroups which are commensurable. Similarly, an element is \emph{commensurable into} a subgroup $H$ if the cyclic subgroup it generates is commensurable into $H$.} elements $a,b\in H$, and $\phi(\{K_1,\ldots,K_n\})=\{K_1,\ldots,K_n\}$ then $\mathcal{R}$ can be chosen so that  $\phi(\mathcal{R})=\mathcal{R}$.
\end{lemma}

\begin{remark}\label{rem: extending slightly sc words}
    Let $G$ and $H,K_1,\ldots,K_n$ be as above, and let $u_1,\ldots,u_m$ be quasigeodesic words in $G$, then if $w_1,\ldots,w_m\in X$ have small enough cancellation and small enough overlap with $K_1,\ldots,K_n$, then so will the words $w_1u_1,\ldots,w_mu_m$.
\end{remark}

\begin{lemma}\label{lem: summary of properties of sc quotients}
    Let $G$ and $K_1,\ldots,K_n$ be as in Lemma \ref{lem: existence of small enough sc words}. Then, for every finite set of words $\mathcal{R}=\{w_1,\ldots,w_m\}$ with small enough cancellation and small enough overlap with $K_1,\ldots,K_n$ the following holds:
    \begin{enumerate}
        \item \label{lem: summary of properties of sc quotients - hyperbolic}The quotient $G/\ngen{\mathcal{R}}$ is torsion-free and hyperbolic.
        \item \label{lem: summary of properties of sc quotients - quasiconvex}For every $1\le i\le n$, the subgroup $K_i$ embeds in $G/\ngen{\mathcal{R}}$ as a quasiconvex subgroup.
        \item \label{lem: summary of properties of sc quotients - non-comm}For every $1\le i,j\le n$, if $K_i$ is non-commensurable into $K_j$ in $G$ then the same holds in $G/\ngen{\mathcal{R}}$.
    \end{enumerate}
\end{lemma}

\subsection{HNN extensions}

The HNN extensions which we use have cyclic edge stabilizers. In this case, one has the following theorem.

\begin{theorem}[Theorem 4 in \cite{kharlampovich1998hyperbolic} or Theorem 1.2 in \cite{kapovich2001combination}] \label{thm: quasiconvexity of cyclic splittings}
    Let $G$ be a hyperbolic group acting on a tree with cyclic edge stabilizers, then the vertex stabilizers of $G$ are quasiconvex. In particular, in hyperbolic HNN extensions with cyclic edge stabilizers, quasiconvex subgroups of vertex groups are quasiconvex in the HNN extension.
\end{theorem}

Since we will need more control over the possible conjugations of elements, we recall the definition of $k$-acylindrical HNN extensions.

\begin{definition}
    Let $k\in \bbN$. An action $G\actson T$ of a group on a tree is $k$-acylindrical if for every $1\ne g\in G$ the fixed-point set of $g$ in $T$ has diameter $\le k$. Equivalently, the pointwise stabilizer in $G$ of a path of length $k$ in $T$ is trivial.

    An HNN extension (and more generally a graph of groups) is $k$-acylindrical if the action on its associated Bass-Serre tree is $k$-acylindrical.
\end{definition}

% \textcolor{red}{N: Sela writes in his paper ``\textit{acylindrical accessibility for groups}'': the HNN extension $G=A*_C = \gen{A,t| tCt\ii=C'}$ is 1-acylindrical if $\forall a\in A, aCa\ii \cap C'=1$. I think it is false. Here is what I believe to be the correct statement:
% The HNN extension $G=A*_C = \gen{A,t| tCt\ii=C'}$ is 2-acylindrical if $\forall a\in A, aCa\ii \cap C'=1$, and $\forall b\in A-C, bCb\ii \cap C=1$. In our case, here is what we need:}

It is easy to verify the following sufficient condition for $2$-acylindricity of a double HNN extension.

\begin{lemma}\label{lem: 2 acylindericity}
    Let $A$ be a group, and $C,C',D,D'$ be distinct subgroups of $A$. Assume that for all $g\in A$, $X\in \{C',D'\}$, and $Y\in \{C,C',D,D'\}$, $gXg\ii \cap Y \ne 1 \implies X=Y, g\in X$. Then, the (double) HNN extension $G=\gen{A,s,t\;|\;C^s = C', D^t=D'}$ $G$ is $2$-acylindrical.\qed
\end{lemma}

Under the condition of $2$-acylindricity it is easy to see the following

\begin{lemma}\label{lem: non-commensurability after HNN}
    Let $A,C,C',D,D',G$ be as in Lemma \ref{lem: 2 acylindericity}, and assume that the edge groups $C,C',D,D'$ are cyclic. Let $U,V$ be two non-commensurable subgroups of $A$. Assume one of the following holds:
    \begin{enumerate}
        \item \label{lem: non-commensurability after HNN - non v cyclic}$U,V$ are not virtually cyclic, or,
        \item \label{lem: non-commensurability after HNN - special case} $U=C$ and $V=D$.
    \end{enumerate}
    Then, $U$ and $V$ are non-commensurable in $G$.
\end{lemma}

\begin{proof}
    If $U$ and $V$ are not virtually cyclic, then the conclusion follows easily from Britton's Lemma and the assumption that the edge groups are cyclic.

    If $U=C$ and $V=D$ then it follows by the Britton's Lemma and the assumption on $C,C',D,D'$ in Lemma \ref{lem: 2 acylindericity} that $C$ and $V$ are not commesurable in $G$.
\end{proof}

% \todo{A lemma that proves \ref{item: y0 not commen into X} and \ref{item: non-commensurability}}

% \begin{proof}
% For property 9:
% Uses the fact that the HNN is over cyclic edge groups, so since $\gen{y,y'}$ and $\gen{x,x'}$ are not cyclic, by Britton's Lemma if an element commensurates $\gen{y,y'}$ into $\gen{x,x'}$ then it must be an element of the vertex group, hence it follows by induction].

% For property 10: any element commensurating $x$ to $y$, must be of length at most 2 by 2-acylindricity, a closer analysis of the possible cases shows that in fact they must be commensurable in the vertex group, hence it is impossible by induction.
% \end{proof}

\section{Proof of Theorem \ref{thm: main}}\label{sec: constructing G}

Theorem \ref{thm: main} follows from the following proposition.

\begin{proposition}\label{prop: main construction}
	There exists a finitely generated non-IG group $G$, an element $x\in G$, and an involution $\phi \in \Aut(G)$ such that for all $g\in G$, $\gen{x^g,\phi(x^g)}=G$.
\end{proposition}

We first prove that it implies Theorem \ref{thm: main}.

\begin{proof}[Proposition \ref{prop: main construction} implies Theorem \ref{thm: main}]
	Let $G$ and $\phi$ be as in Proposition \ref{prop: main construction}. Consider the group $\tilde{G} = G \rtimes \gen{\phi}$. By construction, $G$ contains an index 2 non-IG subgroup. It remains to show that $\tilde{G}$ is FIG. We claim that $\tilde{G}$ is invariably generated by $S=\{x,\phi\}$. That is, $\gen{x^{\tilde{g}},\phi ^ {\tilde{g}'}}=\tilde{G}$ for all $\tilde{g},\tilde{g}'\in \tilde{G}$.
	
	Let $\tilde{H}=\gen{x^{\tilde{g}},\phi ^ {\tilde{g}'}}$. We may assume that $\phi \in \tilde{H}$, by conjugating $H$ by $(\tilde{g}')\ii$ if necessary.
	
	We can write $\tilde{g}=g\phi^{\epsilon}\in\tilde{G}$ where $g\in G$ and $\epsilon \in \{0,1\}$. Since $\phi, x^{g\phi^{\epsilon}}\in \tilde{H}$, it follows that both $x^g$ and $x^{g\phi}=\phi(x^g)$ are in $\tilde{H}$. By the assumption, $G = \gen{ x^g, \phi (x^g)}\subseteq \tilde{H}$, but since also $\phi\in \tilde{H}$ we get that $\tilde{H}=\tilde{G}$.
\end{proof}

\begin{proof}[Proof of Proposition \ref{prop: main construction}]
	The proof is by constructing a group $G$ with the desired properties.
	Let us start with $G(0) = F(x,x',y,y')$, the free group generated by the letters $x,x',y,y'$, and let $\phi\in \Aut(G(0))$ be the involution exchanging $x\leftrightarrow y, x' \leftrightarrow y'$.
	%Denote the subgroups $X=\gen{x,x'}$, $Y=\gen{y,y'}$, and
	Enumerate the elements of $G(0)  = \{g_1,g_2,g_3,\ldots\}$.

	Assume we have constructed a sequence $G(0) \epi G(1) \epi \ldots$ of quotients, $G(n)=G(0)/N_i$ where $N_1\le N_2 \le \ldots$ is an increasing sequence of normal subgroups, and such that the groups $G(n)$ satisfy the following\footnote{we abuse notation and think of elements of $G(0)$ as their images in $G(n)$}:
	\begin{enumerate}[label=(B{\arabic*})]
		\item \label{item: conjugacy complete} The subgroup $\gen{x,x'}$ contains some conjugates of $g_1,\ldots,g_n$.
		\item \label{item: properness} $\gen{x,x'}$ is proper.
		\item \label{item: phi} The automorphism $\phi$ descends to $G(n)$.
		\item \label{item: phi generation} The conjugate $x^{g_n}$ $\phi$-\emph{generates} $G(n)$, i.e, $G(n)=\gen{x^{g_n}, \phi(x^{g_n})}$.
	\end{enumerate}

	Consider the limit $G=\varinjlim G(n) = G(0)/(\bigcup N_n)$. It is a finitely generated group by construction. The subgroup $\gen{x,x'}$ is conjugacy complete by \ref{item: conjugacy complete} and proper by \ref{item: properness},
	%QHere you meant to explain why <x.x'> is proper?
	%(otherwise, $x,x',y,y'\in \gen{x,x'}$ and this can be checked at a finite step of the limit process).
	implying that $G$ is non-IG.
	In addition, $\phi$ is an involution of $G$ by \ref{item: phi}, and for all $g\in G$, $\gen{x^g,\phi(x^g)}=G$ by \ref{item: phi generation}.
	
	To complete the proof of Proposition \ref{prop: main construction} it remains to construct a sequence of quotients as above. To build the sequence $G(n)$ we will use small cancellation, and therefore we would like to assume more on the groups in the process.
	 \begin{enumerate}[label=(B{\arabic*})]
		\setcounter{enumi}{4}
		\item \label{item: hyperbolic} The group $G(n)$ is a torsion-free hyperbolic group.
		\item \label{item: free} $\gen{x,x'}$ is free and quasiconvex.
		\item \label{item: y0 not commen into X} $\gen{x,x'}$ and $\gen{y,y'}$ are not commensurable.
		\item \label{item: non-commensurability} The elements $x,y$ are non-commensurable. In particular, $\gen{x,y}$ is non-elementary.
		\item \label{item: hexagon}
		(\emph{The Hexagon Property}) If $\xi,\xi'\in \gen{x,x'}$ and $z\in G(n)$ satisfy $\xi ^z = \phi((\xi')^z)$ then $\xi'=\xi^{\pm 1}$.
		%\todo{N: I rephrased it so that it matches the later use of the hexagon condition, I also changed the order of the conditions}
		%If an element $\xi\in \gen{x,x'}$ is conjugate to $\eta \in \gen{y,y'}$ by an element of the form $z^{-1}\phi(z)$, then $\phi(\xi)=\eta^{\pm 1}$ (The hexagon property).
		%\item \label{item: some malnormal element} \todo{might be redundant} There exists some $\bar{x}\in \gen{x,x'}$ such that $\bar{x}$ is not commensurable to $\gen{x,x'}$ by an element outside $\gen{x,x'}$.
	\end{enumerate}

	\begin{remark}\label{rem: remark about conditions in theorem B}
	Note the following:
	\begin{itemize}
	    \item \ref{item: y0 not commen into X} implies \ref{item: properness}. In fact, it follows from \ref{item: y0 not commen into X} that $\gen{x,x'}$ has infinite index in $G(n)$.
	    \item \ref{item: hexagon} implies that if $a,b\in G(n)$ are non-commensurable and $\phi(b)=a$ then $\gen{a,b}$ is not commensurable into $\gen{x,x'}$.
	    Otherwise, there exists $z\in G(n)$, $\xi,\xi'\in \gen{x,x'}$ and $N\in \bbN$ such that $a^N = \xi^z$ and $b^N=(\xi')^z$.
	    Applying $\phi$ on the second equation gives $a^N=\phi((\xi')^z)$, from which $\xi^z = \phi (\xi'^z)$ follows. \ref{item: hexagon} then implies that $\xi' = \xi^{\pm 1}$, contradicting the assumption that $a,b$ are non-commensurable.
	\end{itemize}
    \end{remark}
	
	It is easy to verify that $G(0)$ satisfies the above \ref{item: conjugacy complete} - \ref{item: hexagon}. Note that \ref{item: conjugacy complete} and \ref{item: phi generation} are vacuous for $G(0)$.

	Starting with $G(n-1)$ we will build $G(n)$ in a three step process:
	
	\paragraph{Step 1. Conjugating $g_n$ into $\{x,x'\}$ using HNN.}%\todo{N: added these explanations}
	Let $g=g_n$. If $g=1$, set $G'(n)=G''(n)=G(n-1)$ and skip to Step 3.
	%\todo{N: changed it a bit. It used to say if $g\in\{x,x'\}$ skip to Step 2}
	Otherwise, the assumptions of Lemma \ref{lem: existence of small enough sc words} with $H=\gen{x,x'}, K_1=\gen{g}, K_2=\gen{\phi(g)}, K_3=\gen{y,y'}$ are satisfied by \ref{item: hyperbolic}, \ref{item: free}, \ref{item: y0 not commen into X} and \ref{item: non-commensurability}. Therefore, we can find a word $w\in\gen{x,x'}$  such that $w$ satisfies arbitrarily small cancellation in $G(n-1)$, and has arbitrarily small overlap with $\gen{g},\gen{\phi(g)}$ and $\gen{y,y'}$.
    Since $\phi(w)\in\gen{y,y'}$ it follows that $w,\phi(w)$ satisfy arbitrarily small cancellation and small overlap with $\gen{g},\gen{\phi(g)}$.
	
    Let $G'(n)$ be the (double) HNN extension \[G'(n)=\gen{G(n-1),s,t| g^s = w, \phi(g)^t = \phi(w)}.\] and extend $\phi$ by setting it to exchange $s\leftrightarrow t$.
	
	Even though $G'(n)$ is not a quotient of $G(n-1)$ one can make sense of properties \ref{item: conjugacy complete}-\ref{item: hexagon} for $G'(n)$. %\todo{N: added sentence}
	By the induction hypothesis $g_1,\ldots,g_n$ are conjugate into $\gen{x,x'}$ in $G(n-1)$ and therefore also in $G'(n)$; the new HNN relations also conjugate $g=g_n$ to $\gen{x,x'}$, hence $G'(n)$ satisfies \ref{item: conjugacy complete}.
	It is also immediate that $G'(n)$ satisfies \ref{item: properness},\ref{item: phi}.

	Since $w,\phi(w)$ satisfy arbitrarily small cancellation $G'(n)$, by Remark \ref{rem: extending slightly sc words} we see that $g^s=w$ and $\phi(g)^t = \phi(w)$ are also small cancellation relations (in the hyperbolic group $G(n-1)*F(s,t)$). It follows that $w$ can be chosen so that $G'(n)$ satisfies \ref{item: hyperbolic} by Item \ref{lem: summary of properties of sc quotients - hyperbolic} of Lemma \ref{lem: summary of properties of sc quotients}.
	
	Moreover, the groups $C=\gen{g},D=\gen{\phi(g)},C'=\gen{w},D'=\gen{\phi(w)}$ satisfy the conditions of Lemma \ref{lem: 2 acylindericity} as we know that $E(C')=C'$ and $E(D')=D'$ by the "moreover" part of Lemma \ref{lem: sc in F(x,y) implies sc in G}.
	Therefore the HNN extension $G'(n)$ will satisfy \ref{item: free} by Theorem \ref{thm: quasiconvexity of cyclic splittings}.
	It will also satisfy \ref{item: y0 not commen into X} and \ref{item: non-commensurability} by Cases \ref{lem: non-commensurability after HNN - non v cyclic} and \ref{lem: non-commensurability after HNN - special case} of Lemma \ref{lem: non-commensurability after HNN}.
	The proof that the Hexagon Property \ref{item: hexagon} is preserved is slightly more technical and appears in Lemma \ref{lem: hexagon property under HNN}.
	
	Note that at this point $G'(n)$ is not a quotient of $G(n-1)$, and it satisfies all properties except for \ref{item: phi generation}. In the next step, we introduce new relations to $G'(n)$, to make it a quotient of $G(n-1)$.
	
	\paragraph{Step 2. Absorbing $G'(n)$ in a quotient of $G(i-1)$ using small cancellation.}
	%for large $N$, we have that $x^N,y^N$ generate a free \qc~ subgroup (Lemma \ref{lem: large powers gen a free qc group}), and so the assumptions for Lemma \ref{lem: small cancellation w and phi(w) avoiding K} are satisfied.
    As explained in Remark \ref{rem: remark about conditions in theorem B}, it follows from \ref{item: y0 not commen into X} that $\gen{x,x'}$ has infinite index in $G(n-1)$, and both are quasiconvex in $G'(n)$ by \ref{item: free} and Theorem \ref{thm: quasiconvexity of cyclic splittings}. %\todo{it said $G'(i)$ I changed it to $G(i-1)$ because we use words in $H=G(i-1)$ -- which was also incorrect below}
    Using Lemma \ref{lem: non-commensurability after HNN} we see that the conditions of Lemma \ref{lem: existence of small enough sc words} are satisfied for $H=G(n-1),K_1=\gen{x,x'}, K_2=\gen{y,y'}$ in $G'(n)$.
	%$K_3=\gen{x}$ and $K_4=\gen{y}$.
	Hence, by the ``moreover'' part of the lemma, we can find $u\in\gen{x,x',y,y'}$ such that $u,\phi(u)$ have arbitrarily small cancellation in $G'(n)$, and such that $u, \phi(u)$ have arbitrarily small overlap with the subgroups $\gen{x,x'}$ and $\gen{y,y'}$.
	Set
	\[G''(n)=G'(n)/\ngen{s=u,t=\phi(u)}.\]
	
	By the way it is defined the composition $G(n-1) \into G'(n) \epi G''(n)$ is onto. It also follows that $G''(n)$ satisfies \ref{item: conjugacy complete} and \ref{item: phi}.
	By Remark \ref{rem: extending slightly sc words} the relations $s=u$ and $t=\phi(u)$ can be chosen to satisfy arbitrarily small cancellation and small overlap with $\gen{x,x'}$ and $\gen{y,y'}$. Properties \ref{item: hyperbolic}, \ref{item: free}, \ref{item: y0 not commen into X} and \ref{item: non-commensurability} then follow from Lemma \ref{lem: summary of properties of sc quotients}, and the Hexagon Property \ref{item: hexagon} is postponed to Lemma \ref{lem: hexagon property under small cancellation}. As explained in Remark \ref{rem: remark about conditions in theorem B}, \ref{item: properness} follows.

    At this point, $G''(n)$  is a quotient of $G(n-1)$ that satisfies all properties except for \ref{item: phi generation}, which will be taken care of in the last step of the construction.

	\paragraph{Step 3. Forcing $\phi$-generation using small cancellation.} Recall that we denote $g=g_n$. By \ref{item: non-commensurability} $x,y$ are non-commensurable. It follows that so are $x^g$ and $\phi(x^g)=y^{\phi(g)}$. As exaplained in Remark \ref{rem: remark about conditions in theorem B} it follows from Property \ref{item: hexagon} that $\gen{x^g,\phi(x^g)}$ is not commensurable into $\gen{x,x'}$. Using this and \ref{item: free}, we see that $H=\gen{x^g,\phi(x^g)}$ and $K_1=\gen{x,x'}, K_2=\gen{y,y'}$ satisfy the assumptions for the ``moreover'' part of Lemma \ref{lem: existence of small enough sc words}.
	Hence, there exist $v,v'\in \gen{x^g,\phi(x^g)}$ such that $v,v',\phi(v),\phi(v')$ satisfy arbitrarily small cancellation in $G''(n)$ and have arbitrarily small overlap with $\gen{x,x'}$ and $\gen{y,y'}$.

    In order to take care of property \ref{item: phi generation}, we set
    \[G(n) = G''(n)/\ngen{x=v,x'=v',y=\phi(v),y'=\phi(v')}.\]

    We have $G(n-1)\epi G''(n) \epi G(n)$. It follows from the construction that $G(n)$ satisfies  \ref{item: conjugacy complete}, \ref{item: phi} and \ref{item: phi generation}.
    As in Step 2, Properties \ref{item: hyperbolic}, \ref{item: free}, \ref{item: y0 not commen into X} and \ref{item: non-commensurability} follow from Lemma \ref{lem: summary of properties of sc quotients}. The Hexagon Property \ref{item: hexagon} holds by Lemma \ref{lem: hexagon property under small cancellation}, and \ref{item: properness} follows.
\end{proof}

\section{Proof of Theorem \ref{thm: there exists fg IG not FIG group}}\label{sec: constructing fg IG not FIG group}

In the following section we construct a finitely generated IG group that is not FIG, proving Theorem \ref{thm: there exists fg IG not FIG group}.

Let $F=F(a,b)$ be the free group generated by $a,b$, and $F=\{g_1,g_2\ldots\}$ be an enumeration of its elements.
Assume we have found a function $h:F\times F \to F$, elements $\{r_{ij}\}_{i\ge j}\subseteq F$, and a quotient $F\epi G$ that satisfy:
\begin{enumerate}[label=(P{\arabic*})]
    \item \label{p P. condition for IG} For all $s,t,u\in F,$ $\gen{a^s,b^t,h(s,t)^u}=G$.\footnote{As usual we interpret elements of $F$ as their image under the quotient map in $G$}
    \item \label{p P. condition for non-FIG} for all $n\in\bbN$, $\gen{g_1^{r_{n1}},\ldots,g_n^{r_{nn}}}\ne G$.
\end{enumerate}

It is then easy to see that \ref{p P. condition for IG} implies that $G$ is IG, while \ref{p P. condition for non-FIG} implies that it is not FIG.
We therefore wish to find such data.

We first establish some notation.
Set $G(0)=F$. %Similar to the proof of Proposition \ref{prop: main construction}, we use induction to construct a sequence of quotients $G(0) \epi G(1) \epi \ldots$, the function $h$ and elements $\{r_{ij}\}_{j \le i}$. Finally, $G$ will be the limit $G=\varinjlim G(i)$.
Enumerate
\[F\times F = \{(s_1,t_1),(s_2,t_2),\ldots\},\text{~and}\] %and
\[ (F\times F)\times F= \{((s_{j_1},t_{j_1}),u_1),((s_{j_2},t_{j_2}),u_2),\ldots\}.\]
Let $\fresh{\bbN}=\{\fresh{i}\in \bbN  \;|\; j_{\fresh{i}} \notin \{j_1,\ldots,j_{\fresh{i}-1}\}\}$, i.e the set of indices of the enumeration of $(F\times F) \times F$ for which a pair $(s,t)$ is introduced for the first time.
When using the notation $\fresh{i}$, we implicitly assume that the element $\fresh{i}$ is in the set $\fresh{\bbN}$.

Let $n\geq 1$. In the $n$\textsuperscript{th} step of the induction, we will construct:
\begin{itemize}
    \item A group $G(n)$ which is a quotient $G(n-1)\epi G(n)$;%, and thus a quotient of $G(0)=F$;
    \item An image for the pair $(s_{j_n},t_{j_n})$ under $h$, in case this pair has not yet appeared in a previous level. That is, in case $n\in \fresh{\bbN}$.
    \item Elements $r_{nk}\in F$ for all $1\le k \le n$, and a subgroup $K_n:=\gen{g_1^{r_{n1}},\ldots,g_n^{r_{nn}}}$.
    \item Elements $x_{\fresh{i}n}\in F$ for all $1\le \fresh{i} \le n$.
\end{itemize}

Such that the following properties hold in $G(n)$:

\begin{enumerate}[label=(C{\arabic*})]
		\setcounter{enumi}{-1}

    \item \label{p C. G is hyperbolic} $G(n)$ is a torsion-free hyperbolic group.
    \item \label{p C. inductive condition for IG} $\gen{a^{s_{j_n}},b^{t_{j_n}},h(s_{j_n},t_{j_n})^{u_n}} = G(n)$.
    \item \label{p C. a,b non-comm} $a,b$ are non-commensurable. %in $G(n)$.
    \item \label{p C. h(s,t) is non-commensurable into K}
    For all $1\le \fresh{i} \le n$, $h(s_{j_{\fresh{i}}},t_{j_{\fresh{i}}})$
    %$h_{j_{\fresh{i}}}$
    is not commensurable into $K_1,\ldots,K_{\fresh{i}-1}$.% in $G(n)$.
    \item \label{p C. inductive condition for non-FIG}
    The subgroups $K_1,\dots,K_n$ are quasiconvex and free. Since $G$ is torsion-free but not free, it follows from Stallings' Theorem that $K_1,\ldots,K_n$ have infinite index in $G(i)$, and in particular they are proper. % of $G(n)$. %In particular they are proper subgroups of $G(n)$.
    %For all $1\le i \le n$, $K_i:=\gen{g_1^{r_{i1}},\ldots,g_n^{r_{ii}}}$ are quasiconvex and free subgroups of $G(n)$. In particular they are all proper subgroups of $G(n)$.
    \item \label{p C. a,b not commensurable into K}
    For all $1\leq \fresh{i} \le k \leq n$,
    $x_{\fresh{i}k}\in \gen{a^{s_{j_{\fresh{i}}}},b^{t_{j_{\fresh{i}}}}}$ is not commensurable into $K_k$ in $G(n)$.
\end{enumerate}

Finally, we set $G=\varinjlim G(n)$. Notice that property \ref{p C. inductive condition for IG} for $G(n)$ implies that
    $\gen{a^{s_{j_i}},b^{t_{j_i}},h(s_{j_n},t_{j_n})^{u_i}} = G(n)$ for all $i\leq n$, since $G(n)$ is a quotient of $G(i)$. In particular, we get that \ref{p P. condition for IG} holds for $G$. Furthermore, by the definition of the groups $K_i$, Property \ref{p C. inductive condition for non-FIG} implies \ref{p P. condition for non-FIG} for $G$.

It is easy to see that $G(0)=F$ satisfies the above assumptions. Notice however that most conditions are vacuous in this case, as they are defined for $i\geq 1$ only.

We now describe the inductive step. Suppose we have defined the groups $G(0),\ldots,G(n-1)$ with the auxiliary data described above such that they satisfy \ref{p C. G is hyperbolic} - \ref{p C. a,b not commensurable into K}.

\paragraph{Step 1. Defining $h(s_{j_n},t_{j_n})$.} If $n\notin \fresh{\bbN}$ skip this step. Otherwise, $n\in \fresh{\bbN}$ and hence the image of the pair $(s_{j_n},t_{j_n})$ under $h$ was not previously defined. By Lemma \ref{lem: exists h non-comm into K}, there exists an element in $G(n-1)$ that is not commensurable into $K_1,\ldots,K_{n-1}$. Set $h(s_{j_n},t_{j_n})$ to be such an element.

At this point, \ref{p C. h(s,t) is non-commensurable into K} holds also for $i=n$, in $G(n-1)$.

\paragraph{Step 2. Constructing $G(n)$.}
By the induction hypothesis and Step 1, \ref{p C. h(s,t) is non-commensurable into K} for $1\le \fresh{i}\le n$ and \ref{p C. a,b not commensurable into K} for $1\le \fresh{i} \le k<n$ hold in $G(n-1)$. It follows that $\gen{a^{s_{n}},b^{t_{n}},h(s_{n},t_{n})^{u_{n}}}$ contains an element which is not commensurable into $K_1,\ldots,K_{n-1}$.%\todo{isn't \ref{p C. h(s,t) is non-commensurable into K} enough here? as $h$ itself is non-commens into $K_1,\dots,K_{n-1}$}
%\todo{no. h is only good for $K_1,\ldots,K_i$ for the index in which it was defined, beyond this index one has to use C5 to find the non-commensurable element.}

By Lemma \ref{lem: existence of small enough sc words}, there exist words $w_a,w_b\in H=\gen{a^{s_{j_n}},b^{t_{j,n}},h(s_{j_n},t_{j_n})^{u_n}}$ with arbitrarily small cancellation in $G(n-1)$ and arbitrarily small overlap with $K_1,\ldots,K_{n-1}$, $\gen{a},\gen{b}$, $\{\gen{h(s_{j_i},t_{j_i})}\}_{i\leq n}$, %\ldots,\gen{h(s_{j_{n},t_{j_n})}}$
and $\{\gen{x_{ik}}\}_{\fresh{i}\le k \le n-1}$.
Define \[G(n)=G(n-1)/\ngen{w_a=a,w_b=b}.\]
By Item \ref{lem: summary of properties of sc quotients - hyperbolic} of Lemma \ref{lem: summary of properties of sc quotients}, property \ref{p C. G is hyperbolic} persists under small cancellation quotients, and so it holds in $G(n)$.
Moreover, it follows from the new relations that $\gen{a^{s_{j_n}},b^{t_{j,n}},h^{u_n}_{j_n}}=G(n)$, and so \ref{p C. inductive condition for IG} holds for $G(n)$ as well.
Similarly, properties \ref{p C. a,b non-comm} and \ref{p C. h(s,t) is non-commensurable into K} hold in the quotient $G(n)$ by Item  \ref{lem: summary of properties of sc quotients - non-comm} of Lemma \ref{lem: summary of properties of sc quotients} and the induction hypothesis.

Regarding the other two properties:
For all $1\le i \le n-1$, \ref{p C. inductive condition for non-FIG} holds in $G(n)$ by Item \ref{lem: summary of properties of sc quotients - quasiconvex} of Lemma \ref{lem: summary of properties of sc quotients}, since the relations have small overlap with $K_1,\ldots,K_n$.
Similarly, \ref{p C. a,b not commensurable into K} for $1\le \fresh{i} \leq k \le n -1$ holds in $G(n)$ by Item \ref{lem: summary of properties of sc quotients - non-comm} of Lemma \ref{lem: summary of properties of sc quotients}.

It remains to construct $K_n$ and show \ref{p C. inductive condition for non-FIG} for $i=n$, and \ref{p C. a,b not commensurable into K} for $k=n$. This is done in step 3.

\paragraph{Step 3. Constructing $r_{n1},\dots,r_{nn}$ and $x_{\fresh{i}n}$.}
We have seen that \ref{p C. a,b non-comm} holds in $G(n)$, i.e  $a,b$ are non-commensurable in $G(n)$.
Hence for every $1\le\fresh{i}\le n$, %$X_{\fresh{i}} :=
$\gen{a^{s_{j_{\fresh{i}}}},b^{t_{j_{\fresh{i}}}}}$ is non-elementary.
Let
%$Q_{\fresh{i}}\le X_{\fresh{i}}$
$Q_{\fresh{i}}\le \gen{a^{s_{j_{\fresh{i}}}},b^{t_{j_{\fresh{i}}}}}$ be some non-elementary quasiconvex subgroups which exists by Lemma \ref{lem: large powers gen a free qc group}. %(generated by large powers of $a^{s_{j_{\fresh{i}}}},b^{t_{j_{\fresh{i}}}}$).
By Lemma \ref{lem: constructing free groups} below, find $r_{n1},\dots,r_{nn}$ such that $K_n:=\gen{g_1^{r_{n1}},\dots,g_n^{r_{nn}}}$ is quasiconvex and free, and such that for every $1\leq \fresh{i} \leq n$, $Q_{\fresh{i}}$ is not commensurable into $K_n$.
By Lemma \ref{lem: exists h non-comm into K},
for every $1\leq \fresh{i} \leq n$ there exists $x_{\fresh{i}n}\in Q_{\fresh{i}}$ that is not commensurable into $K_n$.

The choice of $r_{n1},\dots,r_{nn}$ ensures \ref{p C. inductive condition for non-FIG} for $i=n$. Lastly, \ref{p C. a,b not commensurable into K} holds for $k=n$ by the construction of $r_{n1},\ldots,r_{nn}$ and $x_{\fresh{i}n}$.

This completes the proof of Theorem \ref{thm: there exists fg IG not FIG group}.\qed

\begin{lemma}\label{lem: constructing free groups}
 Let $G$ be hyperbolic, let $Q_1,\ldots,Q_m\le G$ be some non-elementary quasiconvex subgroups of $G$, and let $g_1,\ldots,g_n$ be infinite order elements of $G$. Then, there exist $r_1,\ldots,r_n\in G$ such that $K=\gen{g_1^{r_1},\ldots,g_n^{r_n}}$ is a quasiconvex free subgroup and $Q_1,\ldots,Q_m$ are not commensurable into $K$.
\end{lemma}

\begin{proof}
Let $\mu = \min\{ \dim _H (\Lambda Q_1), \ldots, \dim_H(\Lambda Q_n)\}$, where $\dim_H(\Lambda Q_i)$ is the visual dimension of the limit set $\Lambda Q_i$ of $Q_i$ in the visual boundary $\partial G$. Paulin \cite{paulin1997critical} shows that $\dim_H(\Lambda Q_i)$ is equal to the critical exponent of the subgroup $Q_i$, and is thus invariant under conjugation. Since $\Lambda Q_i$ is a commensurability invariant, the visual dimension is also a commensurability invariant.

By choosing elements $r_1,\ldots,r_n\in G$ sparse enough, we can make $\dim_H (\Lambda K) < \mu$ where $K=\gen{g_1^{r_1},\ldots,g_n^{r_n}}$ is a quasiconvex free group.
It follows that $Q_1,\ldots,Q_n$ are not commensurable into $K$, as otherwise $\Lambda Q_i \subseteq \Lambda K$ which will contradict the monotonicity of the Hausdorff dimension.
\end{proof}

%\begin{remark}
%This will take care of \ref{p B. h(s,t) is non-commensurable into K} for $1\leq i < n$.
%\end{remark}

\section{Preliminaries on small cancellations with small overlaps}\label{sec: preliminaries on hyp groups and sc}

Let $G$ be generated by a finite set $S$.
Let $W$ be a word over $S$. We write $\|W\|$ to denote the length of $W$ as a word. We use the same notation, $\|p\|$, to denote the length of a path $p$. We often abuse notation and identify a path in the Cayley graph of $G$ with its label.
For an element $g\in G$, we denote by $|g|$ the distance in $\Gamma(G,S)$ between $g$ and $1_G$.
%to denote the length of a shortest path labeled by a word that equals to $g$ in $G$.

\subsection{Small cancellation conditions}
Recall that a set of words $\mathcal{R}$ is called \emph{symmetrized} if it is closed under taking cyclic permutations and inverses.

\begin{definition}[pieces]\label{def: e-piece (modified)}
Let $\mathcal{R}$ and $\mathcal{K}$ be symmetrized sets of words in $S$, and $\epsilon>0$.
Let $U$ be a subword of a word $R\in \mathcal{R}$.
$U$ is called a \emph{($\mathcal{K},\epsilon$)-piece} if there exists a word $R'\in \mathcal{K}$ such that:
\begin{enumerate}
    \item $R= UV$, $R'=U'V'$ as words, for some words $U',V,V'$;
    \item $U'=C UD$ in $G$ for some words $C,D$ in $S$ such that $\max\{\|C\|,\|D\|\}\leq\epsilon$;
    \item $CRC^{-1}\neq R'$ in $G$.
\end{enumerate}
$U$ is called an \emph{$\epsilon'$-piece} if:
\begin{enumerate}
    \item $R=UVU'V'$, for some $U',V,V'$;
    \item $U'=C U^{\pm 1}D$ in $G$ for some words $C,D$ in $S$ such that $\max\{\|C\|,\|D\|\}\leq\epsilon$;
\end{enumerate}
\end{definition}

\begin{remark}
In case $\mathcal{K}=\mathcal{R}$, a $(\mathcal{K},\epsilon)$-piece is simply called an \emph{$\epsilon$-piece}, and this definition coincides with the usual definition found for example in \cite{olshanskii1993residualing, osin2010small}.
\end{remark}

%\textcolor{red}{and here is the old definition:}
%\begin{definition}[$\epsilon$-piece and $\epsilon'$-piece]\label{def: e-piece}
%Let $\mathcal{R}$ be a symmetrized set of words in $S$. For $\epsilon>0$, a subword $U$ of a word $R\in \mathcal{R}$ is called an \emph{$\epsilon$-piece} if there exists a word $R'\in \mathcal{R}$ such that:
%\begin{enumerate}
%    \item $R=UV$, $R'=U'V'$ for some $U',V,V'$;
%    \item $U'=C UD$ in $G$ for some words $C,D$ in $S$ such that $\max\{\|C\|,\|D\|\}\leq\epsilon$;
%    \item $CRC^{-1}\neq R'$ in $G$.
%\end{enumerate}

%Similarly, a subword $U$ of $R\in\mathcal{R}$ is called an \emph{$\epsilon'$-piece} if:
%\begin{enumerate}
%    \item $R=UVU'V'$, for some $U',V,V'$;
%    \item $U'=C U^{\pm 1}D$ in $G$ for some words $C,D$ in $S$ such that $\max\{\|C\|,\|D\|\}\leq\epsilon$;
%\end{enumerate}

%Lastly, a subwords $U$ of $R\in \mathcal{R}$ is a $(K,\epsilon)$-piece, where $K\leq G$ is a subgroup, if there exists a word $R'\in K$ such that
%\begin{enumerate}
%    \item $R=UV$;
%    \item $R'=CUD$ in $G$ for some words $C,D$ in $S$     such that $\max\{\|C\|,\|D\|\}\leq \epsilon$;
%\end{enumerate}
%\end{definition}

\begin{definition}[Small cancellation conditions]\label{def: small cancellation avoiding K}
Let $\mathcal{R}$ and $\mathcal{K}$ be symmetrized sets of words in $G$. We say that $\mathcal{R}$ satisfies the
$\torsionsccK{\mathcal{K}}$ condition for some $\epsilon\geq 0,\mu>0,\lambda\in(0,1],c\geq 0,\rho>0$, if
\begin{enumerate}
    \item \label{item: sc - long words} $\|R\|\geq \rho$ for any $R\in \mathcal{R}$.
    \item \label{item: sc - qg words} any word $R\in \mathcal{R}$ is \qg{$\lambda,c$}, that is, for every subword $V$ of $R$ we have $|V|\ge \lambda\|V\|-c$.
    \item \label{item: sc - no small pieces} for any $(\mathcal{R},\epsilon)$-piece $U$ of any word $R \in \mathcal{R}$, $\max\{\|U\|,\|U'\|\}<\mu\|R\|$.
    \item \label{item: sc - no small K-pieces} for any $(\mathcal{K},\epsilon)$-piece $U$ of any word $R \in \mathcal{R}$, $\max\{\|U\|,\|U'\|\}<\mu\|R\|$.
    \item \label{item: sc - no small epsilon'-pieces} for any $\epsilon'$-piece $U$ of any word $R \in \mathcal{R}$,  $\max\{\|U\|,\|U'\|\}<\mu\|R\|$.
\end{enumerate}
%If all conditions but \ref{item: sc - no small epsilon'-pieces} hold, we say that $\mathcal{R}$ satisfies $\sccK{\mathcal{K}}$.
\end{definition}

\begin{remark}
\begin{enumerate}
    \item %Recall that the \emph{symmetrized closure} of a set of words $E$ is the set of all cyclic permutations of elements of $E$ and $E^{-1}$.
    An arbitrary set of words $E$ is said to satisfy $\torsionscc$ if its symmetrized closure does.
    \item When $\mathcal{K}=\{1\}$, condition \ref{item: sc - no small K-pieces} trivially holds, and the $\torsionsccK{\mathcal{K}}$ conditions coincide with the usual $\torsionscc$ conditions found for example in \cite{olshanskii1993residualing,osin2010small}.
\end{enumerate}

\end{remark}

Instead of keeping track of quantifiers, it would be convenient to use the following.
\begin{definition}
     Let $G,\mathcal{K}$ as in the definitions above. Let $\mathcal{P}$ be some property.

     We say that \emph{there exists a set of words $\mathcal{R}$ satisfying $\mathcal{P}$ in $G$ with arbitrarily small cancellation and arbitrarily small overlap with $\mathcal{K}$} if there exists $\lambda,c$ such that for all $\epsilon,\mu,\rho$ there exists a set $\mathcal{R}$ satisfying $\mathcal{P}$ and the $\torsionsccK{\mathcal{K}}$-condition.

     Similarly, we say that \emph{$\mathcal{P}$ holds for sets of words $\mathcal{R}$ of $G$ with small enough cancellation and small enough overlap with $\mathcal{K}$} if for every $\lambda,c$ there exist $\epsilon,\mu,\rho$ such that $\mathcal{P}$ holds for all $\mathcal{R}$ satisfying the $\torsionsccK{\mathcal{K}}$-condition.
\end{definition}

\begin{remark}
    Suppose $G$ is hyperbolic and $K_1,\ldots,K_n$ are quasiconvex in $G$. Fix some generating sets $S_1,\ldots,S_n,S$ for $K_1,\ldots,K_n,G$ respectively. We assume $S$ contains $S_1,\ldots,S_n$.
    By ``small overlap with $K_1,\ldots,K_n$'' we mean ``small overlap with $\mathcal{K}$'' where $\mathcal{K}=\mathcal{K}_1\cup\ldots\cup\mathcal{K}_n$ and $\mathcal{K}_i$ is the set of all words in $S_i$ which are geodesic in $K_i$.
\end{remark}

%\textcolor{red}{and the old definition:}
%\begin{definition}[Small cancellation condition]\label{def: small cancellation}
%We say that the set $\mathcal{R}$ satisfies the
%$\scc$ condition for some $\epsilon\geq 0,\mu>0,\lambda>0,c\geq 0,\rho>0$, if
%\begin{enumerate}
%    \item \label{item: sc - long words} $\|R\|\geq \rho$ for any $R\in \mathcal{R}$.
%    \item \label{item: sc - qg words} any word $R\in \mathcal{R}$ is \qg{$\lambda,c$}.
%    \item \label{item: sc - no small pieces} for any $\epsilon$-piece $U$ of any word $R \in \mathcal{R}$, the inequality $\max\{\|U\|,\|U'\|\}<\mu\|R\|$ holds.
%\end{enumerate}

%The set $\mathcal{R}$ satisfies the $\torsionscc$ condition if Item \ref{item: sc - no small pieces} also holds for any $\epsilon'$-piece of any word $R\in \mathcal{R}$;

%Lastly, $\mathcal{R}$ satisfies the $C^K(\epsilon,\mu,\lambda,c,\rho)$ condition for some $K\leq G$ if Item \ref{item: sc - no small pieces} is replaced by
%\begin{enumerate}
%\setcounter{enumi}{3}
%    \item for any $(K,\epsilon)$-piece $U$ of any word $R \in \mathcal{R}$, we have $\|U\|<\mu\|R\|$.
%\end{enumerate}

%\end{definition}

\subsection{The Greendlinger Lemma}
%\paragraph{contiguity diagrams}
Let $G=\gen{S|\mathcal{O}}$ be a presentation of $G$,
$\mathcal{R}$ a set of words and $G'= \gen{S|\mathcal{O}\cup \mathcal{R}}$.
Let $\Delta$ be a van Kampen diagram over $G'= \gen{S|\mathcal{O}\cup \mathcal{R}}$ and $q$ a subpath of $\partial \Delta$. Let $\Pi$ be an $\mathcal{R}$-cell of $\Delta$, i.e., a cell whose boundary is labelled by a word in $\mathcal{R}$.
Suppose $\Gamma$ is a subdiagram of $\Delta$, containing no $\mathcal{R}$-cells, and such that $\partial \Gamma=s_1 q_1 s_2 q_2$ where $q_1$ is a subpath of $\partial \Pi$, $q_2$ a subpath of $q$ and $\max \{|s_1|,|s_2| \}\leq \epsilon$ for some $\epsilon >0$.
Then $\Gamma$ is called an \emph{$\epsilon$-contiguity subdiagram of $\Pi$ to $q$},
%and $q_1$ is the \emph{contiguity arc of $\Gamma$ to $q$}.
and the ratio $\|q_1\|/\|\partial \Pi\|$ is called \emph{the contiguity degree} of $\Pi$ to $q$, denoted by $(\Pi, \Gamma, q)$.
%($\Pi,\Gamma, \Sigma$).

%A van Kampen diagram over $G_1= \gen{S|\mathcal{O}\cup \mathcal{R}}$ is said to be reduced if it has minimal number of $\mathcal{R}$-cells among all diagrams over this presentation having the same boundary label.

%\begin{lemma}[geometric reformulation of small cancellation conditions]\label{lem: geometric reformulation of sc}
%Suppose $\mathcal{R}$ satisfies $\scc$ for some values of the parameters. Let $\Delta$ be a reduced diagram over $G_1=\gen{\mathcal{A}|\mathcal{O}\cup \mathcal{R}}$.
%Then for every $\epsilon$-contiguity subdiagram $\Gamma$ of a cell $\Pi$ to another cell $\Sigma$, we have $$\max\{(\Pi,\Gamma,\Sigma),(\Sigma,\Gamma,\Pi)<\mu$$
%\end{lemma}

Let $\Sigma,\Sigma'$ be subdiagrams of $\Delta$ containing no $\mathcal{R}$-cells and such that $\partial \Sigma$ and $\partial \Sigma'$ have the same label.
In this case, replacing $\Sigma$ by $\Sigma'$ will not affect the label of $\partial \Delta$ and the number of $\mathcal{R}$-cells in $\Delta$.
Diagrams over $\gen{S|\mathcal{O}\cup \mathcal{R}}$ that can be obtained from each other by a sequence of such replacements are called
$\mathcal{O}$-equivalent.

The following is an analogue to the well-known Greendlinger's Lemma, proved in
%by Ol'shanskii \cite[Lemma~6.6]{olshanskii1993residualing} and
Osin \cite[Lemma~4.4,~5.1~and~6.3]{osin2010small}.

%\begin{lemma}\label{lem: Greendlinger}
%Let $G=\gen{S|\mathcal{O}}$ be hyperbolic.
%Then for any $\lambda\in (0,1],c\geq 0$ and $\mu\in (0,\frac{1}{16}]$ there exist $\epsilon \geq 0$ and $\rho > 0$ with the following property.
%Let $\mathcal{R}$ be a symmetrized set of words satisfying $\torsionscc$ and
%$\Delta$ a reduced van-Kampen diagram over $\gen{S|\mathcal{O} \cup \mathcal{R}}$ whose boundary is \qg{$\lambda,c$}.
%Assume that $\Delta$ has at least one $\mathcal{R}$-cell.
%Then there exists a diagram $\Delta'$ which is $\mathcal{O}$-equivalent to $\Delta$, an $\mathcal{R}$-cell $\Pi$ in $\Delta'$ and an $\epsilon$-contiguity subdiagram $\Gamma$ of $\Pi$ to $\partial \Delta'$ such that $$(\Pi,\Gamma,\partial \Delta')>1-13\mu.$$
%\end{lemma}

\begin{lemma}\label{lem: Greendlinger}
Let $G=\gen{S|\mathcal{O}}$ be hyperbolic and torsion-free. Then for any $\lambda\in (0,1]$ and $c\geq 0$ there exist $\mu>0$, $\epsilon \geq 0$ and $\rho > 0$ with the following property.
Let $\mathcal{R}$ be a symmetrized set of words satisfying $\torsionscc$ and $\Delta$ a reduced van-Kampen diagram over $\gen{S|\mathcal{O} \cup \mathcal{R}}$ whose boundary is \qg{$\lambda,c$}.
Assume that $\Delta$ has at least one $\mathcal{R}$-cell.
Then there exists a diagram $\Delta'$ which is $\mathcal{O}$-equivalent to $\Delta$, an $\mathcal{R}$-cell $\Pi$ in $\Delta'$ and an $\epsilon$-contiguity subdiagram $\Gamma$ of $\Pi$ to $\partial \Delta'$ such that $$(\Pi,\Gamma,\partial \Delta')>1-13\mu.$$
\end{lemma}
%\todo{shall $\mathcal{K}$ be mentioned here somewhere?}

%\begin{lemma}
%    For any $\lambda\in (0,1]$, $c\geq 0$ there are $\mu>0$, $\epsilon \geq 0$ and $\rho>0$ such that the following condition holds. Suppose that $\mathcal{R}$ is a symmetrized set of words in $A$ satisfying the $\torsionscc$ condition. Then every element of finite order in the group $G_1$ is the image of an element of finite order of $G$.
%\end{lemma}

\section{Existence of small cancellation words}\label{sec: finding sc words}

The goal of this section is to prove Lemma \ref{lem: existence of small enough sc words} which states that there exist words with arbitrarily small cancellation and arbitrarily small overlap with a finite union of quasiconvex subgroups.

\subsection{Quasiconvex subgroups}

We begin by collecting some properties of quasiconvex subgroups in hyperbolic groups.

\begin{lemma}\label{lem: only finitely many K<H up to H-conj}
Let $H,K$ be \qc~ subgroups. The collection $\{H\cap K^g| g\in G\}$ of subgroups of $H$ has finitely many ($H$-)conjugacy classes of subgroup.
\end{lemma}

\begin{proof}
By quasiconvexity of $H$ and $K$, there exists $D$ such that for every conjugate $K^g$ of $K$ that has infinite intersection with $H$, the coset $g\ii K$ must lie within distance $D$ from $H$. Hence, $K^g$ is conjugated to $K$ by some $g'=dh$ where $|d|\leq D$ and $h\in H$.
\end{proof}

\begin{lemma}\label{lem: exists h non-comm into K}
Let $H$ be a non-elementary hyperbolic group, and let $Q_1,\ldots,Q_n$ be infinite index quasiconvex subgroups. Then, there exists $h\in H$ which is not commensurable into $Q_1,\ldots,Q_n$.
\end{lemma}
\begin{proof}
Consider the Gromov boundary $\partial H$ with some metric $d$.
Let $L_i = \Lambda Q_i$ be the limit set of $Q_i$ in $\partial H$. Since $Q_i\le H$ is a quasiconvex infinite index subgroup, $L_i$ is a closed meager subset of $\partial H$.
Fix %some
$0<\epsilon<\diam(H)$.
By Corollary 2.5 of \cite{gitik1998widths}, there are finitely many $H$-translates of $L_1,\ldots,L_n$ with diameter $>\epsilon$.
Let $L$ be the union of all of those translates. $L$ is a closed meager set.

Hence, the set $U=\{(x,y)\in (\partial H)^2\;|\;d(x,y)>\epsilon\} \cap (\partial H - L)^2$ is a non-empty open set. Since the set of pairs of endpoints $\{(h^\infty, h^{-\infty})\;|\;h\in H\}$ is dense in $(\partial H)^2$, we can find an element $h$ such that $(h^\infty,h^{-\infty})\in U$. The element $h$ is not commensurable into $Q_1,\ldots,Q_n$ as otherwise the endpoints $h^\infty, h^{-\infty}$ would be in a translate of $L_1,\ldots,L_n$, contradicting the above.
\end{proof}

\begin{corollary}\label{cor: exists h non-comm into K in G}
Let $G$ be a hyperbolic group, let $H,K_1,\ldots,K_k$ be \qc~ subgroups, and suppose that $H$ is not commensurable into any of the $K_i$.
Then there exists $h\in H$ which is not commensurable into any of the $K_i$.
\end{corollary}

\begin{proof}
By Lemma \ref{lem: only finitely many K<H up to H-conj} the collection $\{H\cap K_j^g | g\in G, 1\le j \le k\}$ is finite up to conjugation in $H$. Let $Q_1,\dots,Q_n$ denote representatives (up to conjugation in $H$) of this collection.
By Lemma \ref{lem: exists h non-comm into K} there exists $h\in H$ which is not conjugate in $H$ to any of $Q_i$. In particular, $h$ does not belong to $K_j^g$ for any $g\in G$ and $1\le j \le k$.
\end{proof}

\subsection{Basic geometry of hyperbolic groups}
%with respect to a fix finitely generating set $S=S^{-1}$ and we denote by $\Gamma=\Gamma(G,S)$ is its Cayley graph.
%The following basic properties of hyperbolic groups will be used to prove Lemma \ref{lem: sc in F(x,y) implies sc in G}.
In this subsection we collect some standard lemmas regarding the geometry of hyperbolic groups. The proofs of the following lemmas can be found in Ol'shanskii \cite{olshanskii1993residualing}. Throughout this subsection $G$ is assumed to be a $\delta$-hyperbolic group.

%\begin{definition} A path $p$ in a metric space is said to be \emph{\qg{$\lambda,c$}} for some $\lambda\geq 1, c \geq 0$, if for every subpath $q$ of $p$, $l(q) \leq \lambda d(q_-, q_+) + c$.
%Any word $w$ over $S$ corresponds to a path $p(w)$ in $\Gamma$ from $1$ to $w$ of length $l(w)$. We say that $w$ is \emph{geodesic} if the path $p(w)$ is geodesic in $\Gamma(G, S)$. Similarly, we say that a word $w$ over $S$ is \emph{\qg{$\lambda,c$}} if the path $p(w)$ is \qg{$\lambda,c$} in $\Gamma(G, S)$.
%\end{definition}

%\begin{lemma}\label{lem: qg is close to a geod}
%Let $\lambda\geq 1, c \geq 0$. There is a constant $\nu=\nu(\delta,\lambda,c)\geq 0$ such that for any \qg{$\lambda,c$}~ path $p$ in $\Gamma$ and a geodesic $q$ with $p_-=q_-$, $p_+ = q_+$, one has $p\subset O_\nu(q)$ and $q \subset O_\nu(p)$.
%\end{lemma}

%The following property of a quasigeodesic rectangle one proves by: dividing the rectangle into two quasigeodesic triangles using the diagonal; approximating the quasigeodesics by geodesics (Lemma \ref{lem: qg is close to a geod}; and finally recalling that in a geodesic rectangle, every edge is contained in the $\delta$-neighborhood of the other edges.

\begin{lemma}[Fellow Traveling]\label{lem: fellow-traveling}
Given $\lambda\in(0,1],c\geq 0$ there exists $\delta'\ge 0$ such that for every $\epsilon\geq 0$, there exists $\epsilon'\ge 0$ %$\delta'=\delta'(\delta,\lambda,c)$
with the following property.
If $ p_1 q_1 p_2 q_2$ is a \qg{$\lambda,c$} rectangle and $\|p_1\|,\|p_2\|\leq\epsilon$, then there exist subpaths $q_i'\subset q_i$ of length
$\|q_i'\|>\|q_i\|-\epsilon'$ such that $q_1'$ and $q_2'$ are of Hausdorff distance at most $\delta'$ from each other.
%Given $\lambda\geq 1$ and $c\geq 0$, let $\nu=\nu(\lambda,c)$ be as in Lemma \ref{lem: qg is close to a geod} and denote by $\delta'=2\delta+2\nu$.
%Then for every \qg{$\lambda,c$} rectangle, every side is contained in the union of the $\delta'$-neighborhoods of the other sides.

%In particular, there exists a constant $\tau=\tau(\delta,\lambda,c)$ such that if a \qg{$\lambda,c$} rectangle $ABCD$ has that $\|B\|>4\lambda \max\{\|A\|,\|C\|\}$ and $\|B\|>\tau$, then there exist subpaths $B_0\subset B,D_0\subset D$, each is contained in the $\delta'$-neighborhood of the other, such that both $B_0$ and $D_0$ have length at least $\frac{1}{4}\|B\|$.
\end{lemma}

\begin{lemma} \label{lem: projection almost preserve order}
Let $\delta'>0$, $\lambda \in (0,1]$ and $c\geq 0$, then there exists $d>0$ with the following property.
Let $p,p'$ be \qg{$\lambda,c$} paths of Hausdorff distance at most $\delta'$ from one another.
%, each contained in the $\delta'$ neighborhood of the other.
Let $q$ be a subpath of $p$ and let $q'_-$ and $q'_+$ be projections of $q_-,q_+$ on $p'$ respectively. Namely, $q'_-$ ($q'_+$) is a nearest point to $q_-$ ($q_+$) in $p'$.
If $\|q\|>d$ then $q'_-$ appears before $q'_+$ in $p'$.

%Then there exists $K$ and $d$, with the following property. Let $p,p'$ be \qg{$\lambda,c$} paths, each contained in the $\delta'$ neighborhood of the other. Let $(p)_p =v_0,v_1,v_2,\dots,v_n=(p)_+$ be points on $p$ (that appear in this order) and $p_i$ the subsegment of $p$ bounded by $v_{i-1}$ and $v_i$ ($i=1,\dots n$). Let $v_i'$ be a points in distance at most $\delta'$ from to $v_i$ in $p'$
%and let $o_i$ denote the geodesic between $v_i$ and $v_i'$.
%If $\|p_i\|>d$ for all $i$, then:
%\begin{enumerate}
%    \item $v_0',\dots,v_n'$ appear in $p'$ in this order. Denote by $p_i'$ the subsegment of $p'$ bounded between $v_{i-1}'$ and $v_i'$ ($i=1,\dots, n$);
    %\item  for all $i$, $ \frac{1}{\lambda}l(p_i)-K \leq l(p_i') \leq \lambda l(p_i)+K$;
    %\item  for all $i$, $p_i'$ contains a subpath $p_i''$ such that $p_i''\subset O_{\epsilon}(p_i)$ and $|l(p_i'')-l(p_i')|<K$
%\end{enumerate}
\end{lemma}

%proof: Olshanskii Lemma 1.9.

%A word in $G$ is called \emph{cyclically minimal} if it is the shorter word among all of its conjugates.
%\begin{lemma}\label{lem: cyclically minimal words are (l,c)-qg}
%    There exist parameters $\lambda\in (1,1],c\geq 0$ such that every cyclically minimal word in $G$ is \qg{$\lambda,c$}.
%\end{lemma}

\begin{lemma}\label{lem: large powers gen a free qc group}
Let $x,y\in G$ be non-commensurable elements.
Then there exists $N>0$ such that $\gen{x^N,y^N}\leq G$ is a free quasiconvex subgroup.
\end{lemma}

A group $H$ is called \emph{elementary} if it is virtually cyclic, i.e, contains a finite index cyclic subgroup. When $G$ is hyperbolic, every infinite order element $g\in G$ is contained in a unique maximal elementary subgroup $E(g)\leq G$, %called the \emph{elementarizer of $g$},
which is given by $E(g)=\{x\in G ~|~ \exists n\neq 0:~ xg^nx^{-1}=g^{\pm n} \}$.
If $G$ is moreover torsion-free, then $E(g)$ is cyclic by Stallings' Theorem.

%\begin{lemma}\label{lem: elementarizer}
%\begin{enumerate}
%    \item $g$ is contained in a unique elementary subgroup $E_G(g)=E(g)\leq G$ called the \emph{elementarizer of $g$}.
%    \item The elementarizer is given by $E(g)=\{x\in G ~|~ \exists n\neq 0:~ xg^nx^{-1}=g^{\pm n} \}$.
%    \item The subgroup $E^+(g)=\{x\in G ~|~ \exists n\neq 0:~ xg^nx^{-1}=g^{n} \}\leq E(g)$ is of index $1$ or $2$ in $E(g)$.
%    \item If for some $x\in G$ and $k,l\in \mathbb{Z}$ we have that $xg^k x^{-1}=g^l$, then $k=\pm l$.
%    \item Suppose $G$ is torsion-free, then $E(g)$ is cyclic. \label{item: cyclic elementarizer}
%\end{enumerate}
%\end{lemma}

%\begin{proof}
%1-4 are in Lemmas 1.6 and 1.7 in O'lshanskii.
%5 is by Stalling's theorem: a torsion-free virtually-cyclic group is cyclic, in particular, elementary groups are such.
%\end{proof}

\begin{lemma}\label{lem: x^m and y^n not close together for large m}
Suppose that $G$ is moreover torsion-free, and let $g,h\in G$ be non-trivial elements. There exist constants $M>0$ and $\theta>0$ such that:

If for some $m\geq M$, $xg^my=h^n$ and $\max\{|x|,|y|\}\leq \theta m$, then
%$xgx^{-1},ygy^{-1}\in E(h)$, and in particular,
$g,h$ are commensurable and $g\neq h^{-1}$.
If moreover $g=h$,
%and $|x|,|y|\leq \theta m$
then $x,y\in E(g)$.
%and $x,y\in E^+(g)$ for $n>0$ and $x,y\in E^-(g)$ for $n\leq 0$
%\textcolor{red}{Lastly if $E(g)$ is cyclic, then $g\neq h^{-1}$}.
\end{lemma}

\begin{lemma}[Corner Trimming]\label{lem: corner trimming}
    %Let $G$ be a hyperbolic group.
    For all $\lambda\in(0,1],c\ge 0$ and $k\in \bbN$ there exist $\delta'\ge 0,\lambda'\in(0,1]$ and $c'\ge 0$ such that if $p_1,\ldots,p_k$ are $(\lambda,c)$-quasigeodesic words, then there exist (possibly empty) words $v_1,\ldots,v_{k-1}$ with $\|v_i\|\le \delta'$ and (possibly empty) subwords $p'_1,\ldots,p'_k$ of $p_1,\ldots,p_k$ respectively, such that \[p_1\ldots p_k=p'_1v_1p'_2v_2\ldots v_{k-1}p'_k\]
    in $G$, and the word on the right hand side is a $(\lambda',c')$-quasigeodesic in $G$.%\todo{G: when $p_1\dots p_k$ is a loop you just cut and avoid the corner between $p_k$ and $p_1$?}
\end{lemma}

\begin{proof}
    The case $k=2$ follows from slimness of quasigeodesic triangles in hyperbolic groups, and for $k>2$ it follows by inductively applying the case $k=2$.
\end{proof}

%\begin{proof}
%O'lshanskii, Lemma 2.1.
%\end{proof}

%\begin{lemma}\label{lem: elementarizer}
%    Let $g\in G$ be of infinite order. Suppose $G$ is torsion-free, then $E(g)$ is cyclic.\todo{N: I wouldn't write it as a separate lemma, but rather mention this in the proof when needed}
%\end{lemma}

%\begin{proof}
%By Stalling's theorem, a torsion-free virtually-cyclic group is cyclic.
%\end{proof}

\subsection{Existence of words with arbitrarily small cancellation}

%Let $\mathcal{E}$ be a symmetrized set of words. Following \cite{osin2010small}, we denote by $\mathcal{LE}$ the set of all subwords of words in $\mathcal{E}$.

Given a set of words $\mathcal{R}\subseteq F(X,Y)$ and words $g,h$ in $S$ we denote by $\mathcal{R}(g,h)$ the symmetrized closure of $\{R(g,h)|R\in\mathcal{R}\}$ where $R(g,h)$ is the word obtained by substituting $g,h$ for $X,Y$.

\begin{lemma} \label{lem: sc in F(x,y) implies sc in G}
Let $G$ be a torsion-free hyperbolic group. Let $a,b\in G$ be infinite order elements in $G$ that are non-commensurable.
Let $\lambda\in (0,1],c\geq 0$, and let $\mathcal{K}$ be a symmetrized set of  \qg{$\lambda_0,c_0$} words, that is closed under taking subwords.
Suppose that $a$ is non-commensurable into $\mathcal{K}$.
There exist $\lambda\in (0,1]$ and $c\geq 0$ such that for any $\epsilon\geq 0, \mu>0,\rho >0$, there are $\mu',\rho', N$ with the following property.

%If $\mathcal{R}\subset F(X,Y)$ satisfies $C(0,\mu',1,0,\rho')$  then $\mathcal{R}(a^N,b^N)$ satisfies $\sccK$ condition in $G$;

If a set of words $\mathcal{R}\subset F(X,Y)$ satisfies $C_1(0,\mu',1,0,\rho')$ in $F(X,Y)$,  then $\mathcal{R}(a^N,b^N)$ satisfies $\torsionsccK{\mathcal{K}}$ condition in $G$;

Moreover, for every $R \in \mathcal{R}(a^N,b^N)$ with small enough cancellation, we have that the elementary group $E(R)=\gen{R}$.
\end{lemma}

Given a word $R(X,Y)$, we denote by $\|R(X,Y)\|_F=\|R(X,Y)\|_F$ the norm of $R$ in the free group $F(X,Y)$ with respect to the generating set $X,X^{-1},Y,Y^{-1}$. For words $g,h$ we denote by $\|R(g,h)\|$ the length of a path labeled by $R(g,h)$ in $G$, with respect to the generating set $S$.

\begin{proof}
Let $\lambda,c,N_0$ be such that for any $N>N_0$, any word in $\gen{a^N,b^N}$ is \qg{$\lambda,c$} (Lemma \ref{lem: large powers gen a free qc group}), and moreover, any word in $\mathcal{K}$ is \qg{$\lambda,c$}. Let $\theta$ and $M$ be as in Lemma \ref{lem: x^m and y^n not close together for large m} for the elements $a$ and $b$.
%(Notice that all of $\lambda,c,N_0,\theta$ and $M$, depend on $a,b$ and $\delta$.)
%Let $\delta'$ be as in Lemma \ref{lem: fellow-traveling},
Let $d$ as in Lemma \ref{lem: projection almost preserve order}, %(depending on $\lambda,c$ and $\delta'$).
and $m>\|a\|,\|b\|$.

Consider the cyclic groups $E(a)$ and $E(b)$. %Since $G$ is torsion-free, Stallings' Theorem implies that they are cyclic.
We denote the elements in $E(a)$ (and $E(b)$) by fractional powers of $a$ (resp. $b$). This notation is justified as the generator of $E(a)$ can be thought of as $a^{\frac{1}{r}}$ for some integer $r$, and similarly for $E(b)$.

Let $\epsilon \geq 0, \mu>0, \rho >0$ be arbitrary. Let $\epsilon',\delta'$ be as in Lemma \ref{lem: fellow-traveling}.
Let $t>0$ be a constant with the following property: if an element $a^{L_1}$ in $E(a)$ has length at most $\delta'$, then $L_1<t$; if $b^{L_2}$ in $E(b)$ has length at most $\delta'$, then $L_2<t$.
Set
 \begin{equation}\label{eq: n and N}
    n>\max\{N_0,%m,
    M,2\frac{\delta'}{\theta},\frac{2m}{\theta},2d,2t\},~N=n^2\tag{$E_N$}
  \end{equation}
%$$n>\max\{N_0,m,M,3\frac{\delta'}{\theta},6d\},~N=n^2 ~~(E_N).$$
Let $k$ be a constant such that for every word $R\in F(X,Y)$,
 \begin{equation}\label{eq: k}
    \|R(a^N,b^N)\|>k\|R(X,Y)\|_F\tag{$E_k$}
  \end{equation}
%$$\|R(a^N,b^N)\|>A\|R(x,y)\|_F~~(E_A)$$
Take $\mu',\rho'$ positive constants such that
 \begin{equation}\label{eq: mu'}
    2\mu'<\frac{\mu k}{mN}\tag{$E_{\mu'}$}
  \end{equation}
%$$2\mu'<\frac{\mu A}{mn},~\text{and}$$
%previous version: $$\mu'< \frac{1}{2m}(\mu-\frac{\delta'}{\rho})$$
and
 \begin{equation}\label{eq: rho'}
    \rho'>\max\{\frac{\rho}{k},\frac{1}{\mu'}(\frac{\epsilon'}{mN}+3) \}\tag{$E_{\rho'}$}
  \end{equation}
%$$\rho'>\max\{\frac{\rho}{A},\frac{1}{\mu'}(\frac{\delta'}{mN}+\frac{2}{n}) \}~~(E_{\rho'})$$
%(In the last equation, we can assume that the RHS is positive, as $\rho$ can be enlarged as much as we want.)
Let $\mathcal{R}\subset F(X,Y)$ be a symmetrized set of words satisfying $C_1(0,\mu',1,0,\rho')$,
we wish to show that $\mathcal{R}(a^N,b^N)$ satisfies $\torsionsccK{\mathcal{K}}$ in $G$.

%Recall Definition \ref{def: small cancellation avoiding K}.
It is direct from the construction that
for every $R\in\mathcal{R}$, $\|R(g,h)\| \ge \rho$ and $R(g,h)$ is a $(\lambda,c)$-quasigeodesic.
It remains to show that $\mathcal{R}(a^N,b^N)$ have small $\epsilon$-pieces, $\epsilon'$-pieces, and $(\mathcal{K},\epsilon)$-pieces.
%here when we wish \|R\|>\rho we consider the actual length of the path, and not the distance between $R_-$ and $R_+$, so the length of a path is the same as if it was in the free group, right?
\paragraph{Small $\epsilon$-pieces.} Let $U$ be a maximal $(\mathcal{R}(a^N,b^N),\epsilon)$-piece of a word $R\in \mathcal{R}(a^N,b^N)$, and denote by $\tilde{R}(X,Y)$ the word in $\mathcal{R}(X,Y)$ such that $R$ is a cyclic permutation of $\tilde{R}(a^N,b^N)$. To show that $U$ is small, assume for contradiction
%and suppose $U$ is of maximal length.
 $$\|U\|\geq \mu \|R\|.$$
%we will show that the corresponding word $R(x,y)\in \mathcal{R}(x,y)$ contains a $0$-piece of length at least $\mu'\|R(x,y)\|_F$, contradicting the assumptions.

Let $U'$ be as in Definition \ref{def: e-piece (modified)}. That is, $U'$ is the initial segment of some word $R'\in \mathcal{R}(a^N,b^N)$, and $U'=CUD$ with $||C||,||D||\leq \epsilon$.
Let $\tilde{R'}\in \mathcal{R}(X,Y)$ be such that $R'$ is a cyclic permutation in $G$ of $\tilde{R'}(a^N,b^N)$.
By Lemma \ref{lem: fellow-traveling} there exists a subpath of $U$ of length at least $\|U\|-\epsilon'$ that is of Hausdorff distance at most $\delta'$ from $U'$. Let $U_0$ be such a subpath, of maximal length. See Figure \ref{fig: sc_words}.
%The choice of $\mu'$ and the fact $\|R\|>\rho$ imply that $ (\mu-\frac{\delta'}{\|R\|})>(\mu-\frac{\delta'}{\rho})>2m\mu'$ and in particular, $$\|U_0\|>\mu \|R\|-\delta'>2m\mu' \|R\| ~~ (*).$$

\begin{figure}
  \centering
  \def\svgwidth{\textwidth}
  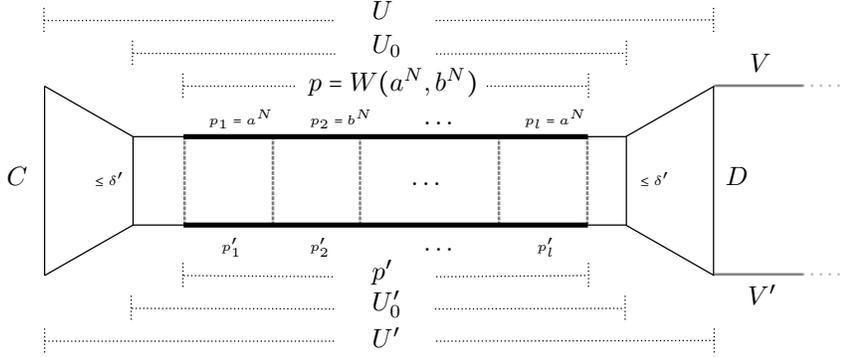
  \label{fig: sc_words}
  \caption{Piece in $\mathcal{R}(a^N,b^N)$}
\end{figure}

Let $W(X,Y)\subset \tilde{R}(X,Y)$ be a maximal word for which $W(a^N,b^N)$ is the label of a subpath $p$ of $U_0$.
Write $p$ as $p=p_1p_2\dots p_l$, according to the letters in $W(X,Y)$. That is, each $p_i$ is labeled $x_i^N$, $x_i \in \{a^{\pm 1},b^{\pm 1}\}$.
%We denote by $Y_i$ the complement letter to $X_i$, $Y_i:=\{a,b\}\setminus X_i$
%Further, part each $p_i$ as $p_i=q_{i,1}q_{i,2}\dots q_{i,n}$ where $q_{i,j}$ is labeled by $X_i^n$.
%Let $p=p_1 p_2\dots p_{r}$ be a maximal subpath of $U_0$ with the following property.\todo{make clear}
%$r=ln$ for some $l\in \mathbb{N}$; Each $p_i$ is labeled by $(X_i)^n$ where $X_i \in \{a,b\}$; For $k=0,\dots ,l-1$, $X_{kn+1}=X_{kn+2}=\dots X_{kn+n}$.
%As $\|(X_i)^n\|<mn$, it follows from $(*)$ that
%$$ln>\frac{2\mu' \|R\|}{n}-2~~ (**).$$
%minus 2 in the RHS here is because $U_0$ may not start and end with a "non-complete block".
As $\|U_0\|\geq \mu \|R\|-\epsilon'$ and $\|x_i^N\|<mN$, we have that

 \begin{equation}\label{eq: l}
    l>\frac{\mu\|R\|-\epsilon'}{mN}-2\tag{$E_l$}
  \end{equation}
%$$ln>\frac{\mu\|R\|-\delta'}{mn}-2~~(*).$$

For each $i$
%(and $j$)
let $p_i'$
%($q_{i,j}'$)
be the projection of $p_i$ on $U'$
%($q_{i,j}$),
%and $p_1'',\dots,p_r''$
%be the projections of $p_1,\dots,p_r$ on $U'$
as explained in Lemma \ref{lem: projection almost preserve order}.
%Suppose without loss of generality that $X_1=a$,
%We now analyze the option for the label of $p_i'$.
Let $y_i$ be an arbitrary letter in $\{a^{\pm 1},b^{\pm 1}\}\setminus \{x_i\}$.
We claim that if $p_i'$ contains a subpath labeled $y_i^L$, then $|L|\leq n$, and this subpath  appears in one of the ends of $p_i'$.
Indeed, suppose $q'\subset p_i'$ was a subpath labeled $y_i^L$ for $|L|>n$. By Lemma \ref{lem: projection almost preserve order}, and since $\frac{n}{2}>d$, there is a subpath $q_0'$ of $q'$ labeled $y_i^{L_0}$, $|L_0|>\frac{n}{2}$ whose projection on $U$ is contained in $p_i$, and therefore labeled by a power of $x_i$.
Since $x_i,y_i$ are either non-commensurable, or inverse to each other, it follows that the rectangle bounded between $q_0'$ and its projection on $U$ contradicts Lemma \ref{lem: x^m and y^n not close together for large m}. Finally, since the label $y_i^L$,$L<n$ cannot be delimited by two appearances of $x_i$, it must be that $q'$ lies in one of the ends of $p_i'$.
It follows that $p_i'$ is labeled by
$y_i^{L_{i,1}} x_i^{S_i} y_i^{L_{i,2}}$,
where $|L_{i,1}|,|L_{i,2}|\leq n$. Observe further that $|S_i-N|\leq 2t<n$. That is since the geodesics connecting the subpath of $p_i'$ labeled by $x_i^{S_i}$ to its projection on $U_0$ are of length at most $\delta'$, and have labels in $E(x_i)$, by Lemma \ref{lem: x^m and y^n not close together for large m}.

We conclude that $p'=p_1'p_2'\dots p_l'$ is labeled

\begin{equation}\label{eq: p'}
    (y_1^{L_{1,1}}x_1^{S_1}y_1^{L_{1,2}})   (x_2^{L_{2,1}}x_2^{S_2}y_2^{L_{2,2}}) \dots
 (y_l^{L_{l,1}}x_l^{S_l}y_l^{L_{l,2}}) %\tag{$E_{p'}$}
\end{equation}
where $|L_{i,1}|,|L_{i,2}|,|S_i-N|<n$ for all $i$.
However, a word over the alphabet $\{a^N,b^N\}$ that is of the above form,
%$$ (y_1^{L_{1,1}}x_1^{S_1}y_1^{L_{1,2}})   (x_2^{L_{2,1}}x_2^{S_2}y_2^{L_{2,2}}) \dots  (x_l^{L_{l,1}}y_l^{S_l}x_l^{L_{l,2}})$$
must simply equal
$ y^{L_1}(x_1^N x_2^N \dots x_l^N)y^{L_2}$, for some $|L_1|,|L_2|<n$.
%Denote by $p''$ the subpath of $p'$ that is labeled by $W(a^N,b^N)=(x_1^N x_2^N \dots x_l^N)$.

In particular, we obtain that $U'$, like $U$, contains a subpath labeled by $(x_1^N x_2^N \dots x_l^N)=W(a^N,b^N)$. It follows that the segment $W(X,Y)$ appears in both $\tilde{R}(X,Y)$ and $\tilde{R}'(X,Y)$. Observe further that \[||W(X,Y)||_F=l>\mu'\|\tilde{R}(X,Y)\|_F.\]
Indeed,
\begin{align*}
    l-1&>\frac{\mu}{mN}\|R\|-(\frac{\epsilon'}{mN}+3)   \\
    &> 2\mu'\|\tilde{R}(X,Y)\|_F-\mu' \|\tilde{R}(X,Y)\|_F \\
    &> \mu'\|\tilde{R}(X,Y)\|_F.
\end{align*}
%$$l>\frac{2\mu' \|R\|-2n}{N}>\frac{\mu' \|R\|}{N}>\mu'\|R(x,y)\|_F$$
First inequality is a rephrasement of Equation \ref{eq: l}; to get the second recall that the $\frac{\mu}{mN}\|R\|$ is at least $2\mu'\|\tilde{R}(x,y)\|_F$ by Equation \ref{eq: k} and Equation \ref{eq: mu'}, while $(\frac{\epsilon'}{mN} + 2)$ is at most $\mu'\|\tilde{R}(x,y)\|_F$ by Equation \ref{eq: rho'}.

As by assumption $\mathcal{R}(X,Y)$ satisfies $C_1(0,\mu',1,0,\rho')$ in $F(X,Y)$, this implies that $\tilde{R}$ and $\tilde{R'}$ are cyclic permutations of one another.
%As they were both defined only up to cyclic permutation, we can assume that $\tilde{R}=\tilde{R'}$.
It follows that also $R$ and $R'$ are cyclic permutations of one another.

%Since $W(X,Y)$, $U_0$ and $U$ were maximal, we have $R\neq R'$, and so $W(X,Y)$ is a $0$-piece in $R(X,Y)$, of length $l$.
%It is left to show that $l>\mu'\|R(X,Y)\|_F$, and indeed,
%$$l-1>\frac{\mu}{mN}\|R(a^N,b^N)\|-(\frac{\epsilon'}{mN}+3)>
%\mu'\|R(X,Y)\|_F-\mu' \|R(X,Y)\|_F>
%\mu'\|R(X,Y)\|_F.$$
%$$l>\frac{2\mu' \|R\|-2n}{N}>\frac{\mu' \|R\|}{N}>\mu'\|R(x,y)\|_F$$
%First inequality is a rephrasement of Equation \ref{eq: l}; to get the second recall that the $\frac{\mu}{mN}\|R(a^N,b^N)\|$ is at least $2\mu'\|R(x,y)\|_F$ by Equation \ref{eq: k} and Equation \ref{eq: mu'}, while $(\frac{\epsilon'}{mN} + 2)$ is at most $\mu'\|R(x,y)\|_F$ by Equation \ref{eq: rho'}.

Denote by $p''$ the subpath of $p'$ labeled by $(x_1^N x_2^N \dots x_l^N)=W(a^N,b^N)$, and part $p''$ as $p''=p_1''p_2''\dots p_l''$ where $p''_i$ is labeled $x_i^N$. We now show that $p=p''$ coincide.
Suppose without loss of generality that $x_1=a$ and let $i$ be the first index for which $x_i=b^{\pm 1}$. Denote the vertices $o=(p_{i-1})_+=(p_i)_-$ and $o''=(p''_{i-1})_+=(p''_i)_-$, and denote the geodesic connecting them by $\gamma$. Notice that $\|\gamma\|\leq \delta'+mn$. Indeed, one can connect them by projecting $o$ on $(p'_{i-1})_+$ and then walk a little (at most $\max\{\|a^n\|,\|b^n\|\}$) along $p'$.
By Lemma \ref{lem: x^m and y^n not close together for large m}, and since $N>\frac{\|\gamma\|}{\theta}$, the label of $\gamma$ must belong to both $E(a)$ and $E(b)$. However, since $a,b$ are non-commensurable, $E(a)\cap E(b)=\{1\}$ and so $\gamma$ is trivial. It follows that $o=o''$. Now, since the labels of $p$ and $p''$ agree, it must be that $p=p''$ fully coincide as paths.

Recall that $R$ and $R'$ are cyclic permutations of one another.
That is: for some word $\Delta$ in $G$, %we have $R=UV$ and $R'=\Delta UV \Delta^{-1}$ in $G$.
one of $\Delta R = R'\Delta $ or $\Delta R' = R\Delta$ holds \emph{as words}. Without loss of generality suppose $\Delta R = R'\Delta $. Since the two copies of $W(a^N,b^N)$ (the one in $R$ and the one in $R'$) are at most $||\gamma||$ from each other, we have that $\Delta$ is very short relative to $R$.

We wish to show that $\Delta=C$, and therefore $CRC^{-1}=R'$, contradicting the assumption that $U$ is an $\epsilon$-piece. To do that, we need to observe first that $U$ and $U'$ coincide not only on the labels, but actually as paths.

Recall that $U$ was assumed to be maximal. As $p=p''$ coincide, and $R$, $R'$ keep agreeing on the labels, we have that  $U\subset R$ and $U'\subset R'$ keep coincide until the beginning of $R$ (on one side) and the end of $R'$ (on the other side). In particular, the start point of $R$ coincides with a point $v$ in $U'\subset R'$, such that $R'$, read from $v$, is identical to $U$ (as words).
It follows that the initial segment of $U'$ ending at $v$ is labeled by $\Delta$. Finally, reading $\Delta^{-1}C$, from the vertex $v$, one arrives back to $v$. Indeed, $\Delta^{-1}$ read from $v$ arrives at $(U')_-$, and $C$ read from there, ends in $v$ again. Then $\Delta=C$ in $G$, as required, and the claim follows.

\paragraph{Small $\epsilon'$-pieces.} The argument $\epsilon'$-pieces follows the same lines as for $\epsilon$-pieces.

\paragraph{Small $(\mathcal{K},\epsilon)$-pieces.} We now show that no small $(\mathcal{K},\epsilon)$-pieces occur.
To do this, suppose now $U$ is a maximal $(\mathcal{K},\epsilon)$-piece of a word $R\in \mathcal{R}(a^N,b^N)$. As before, suppose $\|U\|\geq \mu \|R\|$ and let $U',U_0,W(X,Y),p$ and $p_i,x_i,p_i',1\leq i \leq l$ be as above. Observe that here as well, Equation \ref{eq: l} holds.
Denote by $p'$ be the projection of $p$ on $U'$.

We claim that $x_i\in \{b,b^{-1}\}$ for all $1\leq i \leq l$. Suppose otherwise, that there exists $i$ for which $x_i\in\{a,a^{-1}\}$. Since $p_i'$ is labeled by a word from $\mathcal{K}$, and $a$ is not commensurable into $\mathcal{K}$, the \qg{$\lambda,c$} rectangle bounded between $p_i$ and $p_i'$ would contradict Lemma $\ref{lem: x^m and y^n not close together for large m}$.
It follows that $W(X,Y)=Y^{\pm l}$.
In particular, $Y^{\pm (l-1)}$ is a $0$-piece of $\tilde{R}(X,Y)$, of length $l-1$.
However, the last calculation of the lower bound for $l$, showing $l-1>\mu'\|\tilde{R}(X,Y)\|_F$, holds here as well, contradicting the assumption that $\mathcal{R}(X,Y)$ satisfies $C_1(0,\mu',1,0,\rho')$.

\paragraph{The elementary group $E(R)$ is $\gen{R}$.}
Let now $R\in \mathcal{R}(a^N,b^N)$. %In hyperbolic groups $E(R)$ is virtually cyclic. If $G$ is torsion-free then
Since $G$ is hyperbolic and torsion-free, $E(R)$ is cyclic. In particular, all elements in $E(R)$ commute with $R$.

Let $z\in E(R)$, and consider the quasigeodesic rectangle $t_1u_1=u_2t_2$ representing the relation $z R^B = R^B z$. That is, each of $t_1$ and $t_2$ is labeled by $z$, and each of $u_1$ and $u_2$ by $R^B$.
Recall that $R^B$ is a $(\lambda,c)$-quasigeodesic. By Lemma \ref{lem: fellow-traveling}, there exist $\epsilon',\delta''>0$ depending on $\lambda,c$ and $|v|$ such that $u_i$ have a subpath $u_i'$ of length $\|u_i\|-\epsilon'$ and such that $u_1'$ and $u_2'$ are of Hausdorff distance at most $\delta''$.
Take $B$ large enough, so that $u_1'$ contains a subpath $u_1''$ labeled by $R$.

As in the proof for $\epsilon$-pieces, one shows that the projection of $u_1''$ on $u_2$ is labeled similarly to $u_1''$, and therefore that $u_1''$ and its projection actually coincide in the graph. In particular, since $R$ satisfies small enough cancellation conditions, it must be that the two copies of $R$ have the same `phase'. More precisely, the point $(u_1')_-$ belongs to $u_2$, and the label of $u_2$, read from $(u_1')_-$, starts by $R$ (rather than by a cyclic permutation of it).

Going back to the rectangle $t_1u_1=u_2t_2$, and considering the common point $(u_1')_-=(u_2')_-$, we obtain a triangle $t_1 \bar{u_1}=\bar{u_2}$, where $\bar{u_i}$ is the initial subpath of $u_i$, ending at $(u_1')_-$.
Since both $u_1$ and $u_2$ read from $(u_1')_-$ have labels starting with $R$, it must be that the label of $\bar{u_i}$ is $R^{r_i}$ for some integers $r_i$.
The boundary of the triangle then gives the relation $zR^{r_1}=R^{r_2}$, and it follows that $z\in \gen{R}$.

\end{proof}

As a corollary we can now prove Lemma \ref{lem: existence of small enough sc words}.

%Recall the statement:
%Let $G$ be a torsion-free hyperbolic group, and let $H,K_1,\ldots,K_n$ be quasiconvex subgroups of $G$. If $H$ is non-elementary and non-commensurable into $\mathcal{K} = K_1 \cup \ldots \cup K_n$.
%Then for every $m$ there exists a subset of $m$ words $\mathcal{R}=\{w_1,\ldots,w_m\}\subseteq H$ with arbitrarily small cancellation and arbitrarily small overlap with $\mathcal{K}$.

%If moreover $G$ has an involution $\phi$ which exchanges two non-commensurable elements $a,b\in H$, and $\phi(\mathcal{K})=\mathcal{K}$ then $\mathcal{R}$ can be chosen so that  $\phi(\mathcal{R})=\mathcal{R}$.

\begin{proof}[proof of Lemma \ref{lem: existence of small enough sc words}.]
Say we are given $H,K_1,\dots,K_n$ as in the statement of the Lemma. Since $H$ is non-elementary, we can find $a,b\in H$ that are non-commensurable. Suppose without loss of generality that the generators of each of $K_1,\dots,K_n$ belong to $S$.
Since $K_1,\dots,K_n$ are quasiconvex, the set $\mathcal{K}$ of all elements in $K_1 \cup \dots \cup K_n$ is closed under taking subwords, and all words in $\mathcal{K}$ are \qg{$\lambda_0,c_0$} with respect to some uniform $\lambda_0\in(0,1],c_0\geq 0$.

The first part of \ref{lem: existence of small enough sc words} then follows immediately from Lemma \ref{lem: sc in F(x,y) implies sc in G}. Indeed, given parameters $(\epsilon,\mu,\lambda,c,\rho)$, it is enough to construct arbitrarily large sets of words %$\mathcal{R}=\{W_1,\ldots,W_m\}\subset F(X,Y)$
satisfying $C_1(0,\mu',1,0,\rho')$ in the free group $F(X,Y)$. Such sets are easy to construct. For example, take $N>\max\{\rho',\frac{3}{\mu'}\}$, and for $1\leq i\leq m$ set %$\mathcal{R}:=\{W_1,\dots,W_m,W_1',\dots,W_m'\}$, where $W_i,W_i'$ are defined by
$$ W_{i}=X^{iN} Y X^{iN+1} Y X^{iN+2} Y\dots X^{iN+N} Y,\text{~and}$$
$$ W_{i}'=Y^{iN} X Y^{iN+1} X Y^{iN+2} X\dots Y^{iN+N} X.$$

For the "moreover" part, suppose $\phi$ is an involution of $G$ exchanging two non-commensurable elements $a,b\in H$, and suppose further that $\mathcal{K}=\phi(\mathcal{K})$.
%We wish to choose a family of words $\mathcal{R}$ over the alphabet $a,b$ so that $\phi(\mathcal{R})=\mathcal{R}$.
%Let $H\le G$ be a quasiconvex free subgroup. Let $\phi$ be an involution of $G$. Suppose $H$ is generated by $a,b$ which are non-commensurable in $G$ and $\phi$ exchanges $a\leftrightarrow b$.
%Let $K$ be a quasiconvex subgroup of $G$, such that $H$ is not commensurable into $K$.
%Let $K_1,\ldots,K_n$ be quasiconvex subgroups of $G$, such that $H$ is not commensurable into $K_i$
It is enough to find elements $a',b' \in H$ non-commensurable in $G$, such that $\phi$ exchanges $a'\leftrightarrow b'$ and such that $a'$ is non-commensurable into $\mathcal{K}$.
Indeed, given such elements, one can then apply Lemma \ref{lem: sc in F(x,y) implies sc in G} with $a',b'$, and take the words $W_1,\dots,W_m,W_1',\dots,W_m'$ as suggested above.

We will now find such elements.
Let $h\in H$ be an element not commensurable into $\mathcal{K}'=\mathcal{K}\cup\gen{a}\cup \gen{b}$. %Since $\mathcal{K}=\phi(\mathcal{K}$, also $\phi(g)$ is non-commensurable into $\mathcal{K}$.
For large enough integers $s,S$, the elements
$a'= (a^s h^s)^S$ and
$b'= (b^s \phi(h)^s)^S$
satisfy the requirements.
Indeed, %$\phi(a')=b'$ be construction.
suppose that for some integer $l$ and $g\in G$ we had that  $g^{-1}a'^lg=U$ is either a power of $b'$ or a word in $\mathcal{K}$. We may assume that $a'^l$ is much longer than $g$, by replacing $l$ by a large multiple.
%suppose $g^{-1}a'^lg=k\in \mathcal{K}$ for some integer $l$ and $g\in G$.
By Lemma \ref{lem: fellow-traveling} there exists a major part of $a'^{l}$ that is contained in a small neighborhood of $U$. In particular, by largeness of $S$, this major part must contain a subpath labeled by $a^s h^s$. However, for $s$ large enough, this is impossible by Lemma \ref{lem: x^m and y^n not close together for large m}, as $a$ is non-commensurable with $b$ and $\phi(h)$, and $h$ is non-commensurable into $\mathcal{K}$.
\end{proof}

\subsection{Properties of small cancellation quotients}

In this subsection we prove Lemma \ref{lem: summary of properties of sc quotients} which listed three properties of small cancellation quotients.

\begin{proof}[Proof of Lemma \ref{lem: summary of properties of sc quotients}]

Let $G$ be hyperbolic, let $K_1,\ldots,K_n\le G$ be quasiconvex subgroups, let $\mathcal{R}\subseteq G$ be a symmetrized finite collection of words satisfying small enough cancellation and small enough overlap with $K_1,\ldots,K_n$.
Let $G'=G/\ngen{\mathcal{R}}$.

\paragraph{1. $G'$ is a torsion-free hyperbolic group.} This is proved in \cite{olshanskii1993residualing}.

\paragraph{2. $K_i$ are embedded in $G'$ as quasiconvex subgroups.}
We will show that for every $\lambda\in (0,1],c\geq 0$ there exist $\epsilon\geq0,\mu>0,\rho>0, \lambda'\in (0,1], c'\geq 0$ such that if $\mathcal{R}$ satisfies $\torsionsccK{\mathcal{K}}$-condition, then every $k\in K_i$ that is $(\lambda,c)$-quasigeodesic in $G$ is $(\lambda',c')$-quasigeodesic in $G'$.

It suffices to prove that $|k|>\lambda'\|k\|-c'$ for some $\lambda'$ and $c'$ that are independent of $k$.
Let $k=g$ for some word $g$ which is a geodesic in $G'$.

We would like to apply the Greendlinger Lemma to the relation $k=g$, however, the word $kg\ii$ might not be a qausi-geodesic. Since $k$ and $g$ are quasigeodesics, the only problem that could happen is that there might be a `quasi-backtracking' between $k$ and $g\ii$.
To fix this, we perform corner trimming. By Lemma \ref{lem: corner trimming}, there exist $\delta'>0,0<\lambda'<\lambda,c'>c$ (in what follows, all constants depend only on other constants and never on specific paths) and words $k',v',g'$ such that $k'v'(g')\ii$ is a $(\lambda',c')$-quasigeodesic, $k',g'$ are subwords of $k,g$ respectively, $\|v'\|<\delta'$, and $kg\ii = k'v'(g')\ii$.

If $\|g\|<\lambda'\|k\|-c'$ then $k\ne g$ in $G$. Therefore any van-Kampen diagram of the relation $k=g$ must contain an $\mathcal{R}$-cell.
Assume that $\mu$ is small enough, and $\epsilon,\rho$ are large enough (to be determined later) such that the conclusion of Greendlinger's Lemma holds for $(\lambda',c')$-quasigeodesics. Assume $\mathcal{R}$ satisfies $\torsionsccK{\mathcal{K}}$-condition, then in some van-Kampen diagram $\Delta$ for the relation $k=g$ in $G'$ there exists an $\mathcal{R}$-cell $\Pi$ and an $\epsilon$-contiguity subdiagram $\Gamma$ of $\Pi$ to $\Delta$ such that $(\Pi,\Gamma,\partial\Delta)>1-13\mu$.%\todo{G: in this paragraph the explanation that there is an $\mathcal{R}$-cell in the relation $g=k$ repeats twice?}\todo{N: The first is general, and the second is about Greendlinger's lemma}

Let $\partial \Gamma = s_1 r' s_2 q'$ where $|s_1|,|s_2|\le \epsilon$ and $r'$ is a subpath of $r:=\partial \Pi$ and $q'$ is a subpath of $q=k'v'(g')\ii=\partial \Delta$. We know that $\|r'\|>(1-13\mu)\|r\|$.
Applying Lemma \ref{lem: fellow-traveling} to the quasigeodesic rectangle $\partial \Gamma$,
there exists $\epsilon'$ and subpaths $r'',q''$ of $r',q'$ of lengths $\|r''\|>\|r'\|-\epsilon'$ and $\|q''\|>\|q'\|-\epsilon'$ which are at Hausdorff distance $\delta''$ apart in $G$.
Let $k''=q''\cap k, g''=q''\cap g\ii$  be (possibly empty) subpaths of $k,g$ respectively. Let $r_k,r_g$ be the subpaths of $r$ which are at distance $\delta''$ from $k'',g''$ respectively, and $r''=r_k r_g$.
Combining the above we get, \[\|r_k\|+\|r_g\| = \|r''\| > \|r'\| -\epsilon' > (1-13\mu) \|r\| - \epsilon'>(1-14\mu)\|r\|\]
where the last inequality follows if $\rho$ is large enough.

Let us choose $\mu$ small enough and $\rho$ large enough so that
\begin{equation}\label{eq: inequality for quasiconvexity}
    15\mu\|r\|+2\delta' < \lambda' ((1-15\mu) \|r\|) -c' -2\delta'.
\end{equation}

We divide into two cases:

Case 1. $\|r_k\|> \mu \|r\|$. In this case, we get a contradiction to the small overlap condition with $K_i$.

Case 2. $\|r_g\| > (1-15\mu)\|r\|$. In this case,
let $t_1,t_2$ be paths of length $\le\delta'$ such that $g''=t_1r_g\ii t_2$. Let $r_c$ be the subpath of $r$ which is complementary to $r_g$, i.e $r$ is a cyclic conjugate of $r_g\ii r_c$. Then $g''=t_1 r_c t_2$ in $G$. But
\begin{align*}
    \|t_1r_ct_2\| &\le \|t_1\| + \|r_c\| + \|t_2\| \\
    &\le 15\mu\|r\|+2\delta'\\
    &< \lambda' ((1-15\mu) \|r\|) -c' -2\delta'\\
    &\le \lambda' \|r_g\| -c' -2\delta'\\
    &\le |r_g| - \|t_1\| - \|t_2\| \le\|g''\|
\end{align*}
where the third inequality is by \eqref{eq: inequality for quasiconvexity} and the fifth inequality is by $(\lambda',c')$-quasiconvexity of $r_g$.  This contradicts the assumption that $g$ is a geodesic, as $t_1r_ct_2$ is a shortcut of a subpath of $g$.

It follows from the above that $K_i$ embeds in $G'$. However, one can also easily prove it directly. Assume $k$ is a quasigeodesic word in $K_i$ such that $k\ne 1\in G$ but $k=1\in G'$. Then by Greendlinger's Lemma, a relation $r\in \mathcal{R}$ must have large contiguity degree with $k$ contradicting the small overlap of $\mathcal{R}$ with $\mathcal{K}$.

% We show that if $k\in K_i$ is a $(\lambda,c)$ quasigeodesic in $G$, then it is a geodesic in $G'$. Otherwise there exists a word $g$ which is a geodesic in $G'$ with $k=g \in G'$, and $g\le \lamb$ is shorter than $k$. Therefore, a van-Kampen diagram for the relation $k=g$ must contain cells which are labelled in $\mathcal{R}$.
% Since $\mathcal{R}$ satisfies \torsionsccK{\mathcal{K}}-condition with some $\lambda,c$ and $\mu,\rho$ to be determined later,
% by Greendlinger's Lemma, there must be a relation $r\in \mathcal{R}$ with $r=r'r''$ and such that $r'$ is $\epsilon$-contiguous to a subword $q'$ (i.e, $Ur'=q'V$ for some $U,V$ of length at most $\epsilon$) of the boundary relation $kg\ii$, and satisfies $\|r'\|>(1-13\mu) \|r\|$.
% Let $q'=Q_1Q_2$ be such that $Q_1$ is a subpath of $k$ and $Q_2$ is a subpath of $g\ii$.

% Case 1. $\|Q_2\|$
% then it is geodesic (also in $G$).
% Since $r$ is a $(\lambda,c)$-quasigeodesic the subpath \[\|q'\|
% \ge\lambda(1-13\mu)\|r\|-2\epsilon.\]
% By choosing $\mu$ small enough and $\rho$ large enough (recall that $\|r\| > \rho$) we can ensure that $13\mu\|r\| + 2\epsilon <\lambda(1-13\mu)\|r\|-2\epsilon.$
% And it follows that replacing the subpath $q'$ of $g\ii$ by $Ur'' \ii V\ii$ produces a shorter path which equals $g\ii$. Contradicting that $g$ is a geodesic.
% Otherwise, the $\epsilon$-contiguity must have at least $ $ be along the word $k$, contradicting the fact that elements of $\mathcal{R}$ have small overlap with $K_i$.

\paragraph{3. If $K_i$ is non-commensurable into $K_j$ in $G$ then the same holds in $G'$. }
By Corollary \ref{cor: exists h non-comm into K in G} there exists $h\in K_i$ non-commensurable into $K_j$. We will outline the proof that the same holds in $G'$.

Assume that $h$ is commensurable to $K_j$ in $G'$. Then, there exists $g\in G$, which we may assume to be a geodesic in $G'$, such that $gh^ng\ii = k$ for some $n\in \bbN$ and $k\in K_j$.
Without loss of generality we may assume that $h\in K_j$ is cyclically quasigeodesic, and that $n$ and $\|k\|$ are much larger than $\|g\|$.

As in the proof of Item \ref{lem: summary of properties of sc quotients - quasiconvex} above, we wish to apply the Greendlinger Lemma, and so one has to trim the backtracking corners of the path $gh^ng\ii k\ii$. Since this relation does not occur in $G$, by the Greendlinger Lemma, there must be an $\epsilon$-contiguity between a relation $r\in\mathcal{R}$ and the trimmed path of $gh^ng\ii k\ii$. Since $r$ has small overlap with $K_i$ and $K_j$, the contiguity cannot have a long overlap with $h^n$ nor with $k$, as in Case 1 of the proof of Item \ref{lem: summary of properties of sc quotients - quasiconvex}. Since $g$ is geodesic, the $\epsilon$-contiguity cannot have too long of an overlap with $g$, as otherwise one would be able to shortcut as in Case 2 of the proof of Item \ref{lem: summary of properties of sc quotients - quasiconvex}.
\end{proof}

\section{The Hexagon Property}\label{sec: hexagon}

Let $G$ be a group with an involution $\phi$, let $X\le G$ be a subgroup. Recall that $G$ has the \emph{hexagon property} with respect to $X,\phi$ if for all $\xi,\xi'\in X$ and $z\in G$: $\xi ^z = \phi((\xi')^z)$ implies $\xi'=\xi^{\pm 1}$.

\subsection{Hexagon condition for HNN extensions}

\begin{lemma}\label{lem: hexagon property under HNN}
    Let $A$ be a group with an involution $\phi$, $X\le A$ a subgroup. Let $C \le X$ and $C'\le A$ such that $C,C',D=\phi(C),D'=\phi(C')$ satisfy the conditions of Lemma \ref{lem: 2 acylindericity}. Set $G=\gen{A,s,t \;|\; C^s=C',D^t=D}$. Extend $\phi$ to an involution of $G$ by setting $\phi(s)=t$. If $A$ satisfies the hexagon property with respect to $X,\phi$, then so does $G$.
\end{lemma}

\begin{proof}
    Assume $\xi ^z = \phi({\xi'}^z)$, for some $\xi,\xi'\in X$ and $z\in G$.

    Write $z$ in normal form as $z=a_0 x_1 a_1 \ldots x_n a_n\in G$, where $a_i\in A, x_i \in \{s,s\ii,t,t\ii\}$. Without loss of generality, assume that $z$ has the minimal $n$ among all that satisfy $\xi ^z = \phi({\xi'}^z)$.

    By the assumption on $A$, $z\notin A$. Hence, $n\ge 1$. The word $z \phi (z)\ii$ is reduced in the HNN extension.
    By Lemma \ref{lem: 2 acylindericity}, the extension $G$ is 2-acylindrical. It follows that $n\le 1$.

    Write $z=axb$ where $a,b\in A, x\in \{s,s\ii,t,t\ii\}$.
    The relation $\xi ^z = \phi({\xi'}^{z})$ becomes
    \[ b\ii x\ii a\ii \xi  a x b \quad \phi(b\ii x\ii a\ii {\xi'}\ii a x b)  =1 .\]

    By symmetry, there are two cases to consider:

    Case 1: $x=s\ii$.
    Here the relation becomes
    \[ b\ii  \;\; \overbrace{ s \;\; \underbrace{a\ii \xi  a}_{\in A} \;\; s\ii}^{\heartsuit} \;\; \underbrace{b  \phi(b)\ii}_{\in A} \;\; \overbrace{t \;\; \underbrace{\phi(a \ii{\xi'} \ii a)}_{\in A} \;\; t\ii }^{\heartsuit}\;\;\phi(b)  =1 .\]

    By Britton's Lemma, the word must be non-reduced at both expressions marked with $\heartsuit$.

    After reducing and rearranging we get $(c)^b =   (d) ^{\phi(b)}$
    where $c=s a\ii \xi  a s\ii \in C$ and $d=t \phi(a\ii \xi' a) t\ii \in D=\phi(C)$. Since $c\in C\le X$ and $d=\phi(c')$ for some $c'\in C \le X$ we can apply the hexagon condition of $A$ to deduce that $c'=c^{\pm 1}$. Tracing back the definition of $c,c'$, it follows that $\xi'= \xi^{\pm 1}$, as desired.

    Case 2: $x=s$.
    Applying the same argument we get
    $(c')^b =   (d') ^{\phi(b)}$
    for some $c'\in C', d'\in D'$.
    However, this contradicts the assumption that $g C' g\ii \cap D'=1$ for all $g\in A$.
\end{proof}

\subsection{Hexagon property for small cancellation quotients}

\begin{lemma}\label{lem: hexagon property under small cancellation}
    Let $G$ be a torsion-free hyperbolic group with an involution $\phi$, let $X\le G$ be a quasiconvex subgroup. For all $\mathcal{R}$ such that $\phi(\mathcal{R})=\mathcal{R}$ with small enough cancellation and small enough overlap with $X$, if $G$ has the hexagon property with respect to $X,\phi$ then so does $G/\ngen{\mathcal{R}}$.
\end{lemma}
%\textcolor{red}{maybe rephrase: Suppose $G$ satisfies the hexagon property with respect to $X,\phi$. Then for any $\mathcal{R}$ such that $\phi(\mathcal{R})=\mathcal{R}$ with small enough cancellation and small enough overlap with $X$, the quotient  $G/\ngen{\mathcal{R}}$ satisfies the hexagon property as well.}

\begin{proof}
Assume for contradiction that there exist $\xi,\xi'\in X,z\in G$ such that \[\xi ^z = \phi((\xi')^z) \in G/\ngen{\mathcal{R}}\]
but $\xi'\ne\xi^\pm 1$.
Let us assume that $\xi,\xi'$ are $(\lambda,c)$-quasigeodesics in $G$, and that $z$ is a geodesic in $G/\ngen{\mathcal{R}}$.
The word $q:=z\ii \xi z \phi(z)\ii \phi(\xi')\ii \phi(z)$ is trivial in $G'$ but is not trivial in $G$ since $G$ is assumed to satisfy the hexagon property.
We would like to apply Greendlinger's Lemma to the path $q$.
However, even though the path $q$ is a concatenation of 6 quasigeodesic paths in $G$, it might not be a quasigeodesic because of ``backtracking''.
However, one can fix this by trimming the backtracking corners as described in Lemma \ref{lem: corner trimming}. There exist (possibly empty) subwords $z_1,z_2,z_3,z_4$ of $z$ and subwords $\eta,\eta'$ of $\xi,\xi'$ respectively, and words $v_1,\ldots,v_6$ of length $\le \delta'$  such that the path
\[p:= z_1\ii v_1 \eta v_2 z_2 v_3 \phi(z_3)\ii v_4 \phi(\eta')\ii v_5 \phi(z_4) v_6\]
is a conjugate of $q$ in $G$, and the path $p$ is a $(\lambda',c')$-quasigeodesic, where $\delta',\lambda',c'$ depend only on $\lambda,c$ and $G$. See Figure \ref{fig: hexagon}. 
Moreover, by symmetry of $z\phi(z)\ii$ we may assume that $z_2$ and $z_3$ end at the same place in $z$ (i.e, $z=z' z_2 u = z'' z_3 u$ as words, for some $z',z',u$). A similar statement holds for $z_4, z_1$. By replacing $\xi,\xi'$ with large enough powers, we may assume that $\eta$ and $\eta'$ are arbitrarily long, and in particular non-empty.

Since $p$ and $q$ are conjugates, we have that $p=1\in G/\ngen{\mathcal{R}}$ while $p\ne 1$ in $G$.
By Greendlinger's Lemma there exists a cell labeled $r\in\mathcal{R}$ with contiguity degree $>(1-13\mu)$ assuming $\mathcal{R}$ satisfies small enough cancellation.
Let us denote by $r',p'$ the subwords of $r,p$ respectively which label the opposite sides of the contiguity subdiagram.
As in the proof of Item \ref{lem: summary of properties of sc quotients - quasiconvex} of Lemma \ref{lem: summary of properties of sc quotients} let $r'',p''$ be the $\delta''$-fellow-travelling subpaths of $r',p'$ of length $\|r''\| > \|r'\| - \epsilon', \|p''\| > \|p'\| - \epsilon'$ provided by Lemma \ref{lem: fellow-traveling}, and let $r''=r_{z1}r_{\eta}\ldots r_{z4}$, where $r_{z1}, r_{\eta}\ldots, r_{z4}$ are the (possibly empty) subwords of $r''$ which correspond to the paths that $\delta''$-fellow-travel with $z_1\ii,\eta,\ldots,\phi(z_4)$ respectively. Since $\|r''\|>\|r'\|-\epsilon'$,
\[ \| r_{z1}\| + \|r_{\eta}\| + \ldots + \|r_{z4}\| > (1-13\mu)\|r\| - \epsilon'=:\omega\]

We now divide into cases:

\begin{figure}
  \centering
  \def\svgwidth{\textwidth}
  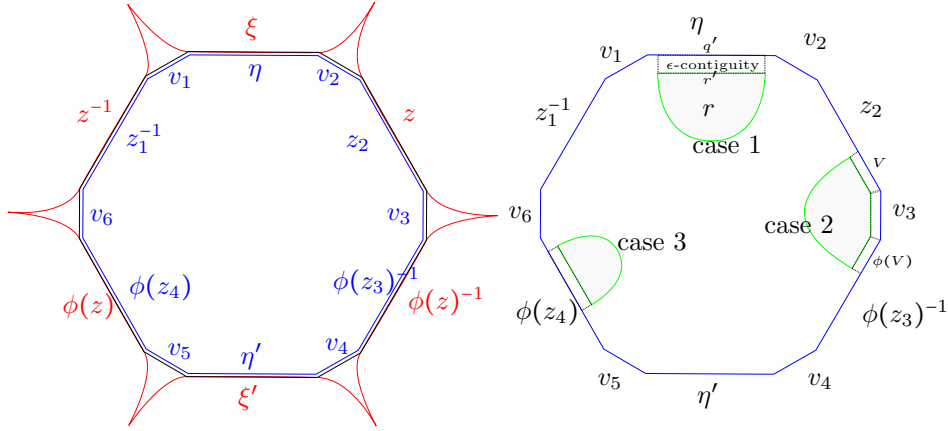
  \label{fig: hexagon}
  \caption{The trimmed hexagon, and the 3 cases of the contiguous cell in the proof of Lemma \ref{lem: hexagon property under small cancellation}}
\end{figure}

Case 1. $\|r_\eta\| > \mu\| r\|=:\omega_1$ or $\|r_{\eta'}\| > \mu \| r\|$. This is impossible when $\mathcal{R}$ has small enough overlap with $X$ since $\eta,\eta'\in X$.

Case 2. $\min \{ \|r_{z2} \|, \| r_{z3}\|\} > \lambda\ii(\mu\|r\|+2\delta''+c)+2\delta''=:\omega_2$ and the path $p''$ contains $v_3$.
In this case, let $p''_{z2} = p''\cap z_2$, $p''_{v3} = p''\cap v_3$ and $p''_{z3} = p''\cap \phi(z_3)\ii$.
Since $\|r\|$ is a $(\lambda,c)$-quasigeodesic and $z_2$ and $\phi(z_3)\ii$ are geodesics we get that $\min\{\|p''_{z2}\|,\|p''_{z3}\|\} > \mu \|r\|+2\delta''$.
Recall that $z_2$ and $z_3$ end at the same place in $z$, thus there is a subword $V$ of $z$ of length $\|V\|>\mu\|r\|+2\delta''$ such that $V$ is in $p''_{z2}$ and $\phi(V)\ii$ is in $p''_{z3}$.
Let $U$ and $U'$ be the subwords of $r$ that  $\delta''$-fellow-travel with $V$ and $\phi(V)$. $\|U\|\ge \|V\| - 2\delta'' > \mu \|r\|$ and similarly $\|U'\|  > \mu\|r\|$.
Since $\phi(r)\in \mathcal{R}$ we get that $r$ has a $2\delta''$-piece (and hence an $\epsilon$-piece) with $\phi(r)$ of length $>\mu \|r\|$ which is impossible if $\mathcal{R}$ has $\torsionsccK{X}$.

% end at the same place in $z$, it follows that there exist words $U,V,W$ such that $z=UVW$ (as words) and

% $p''$ contains both $V$ and $\phi(V)\ii$.
% Thus, $r$ has an $\epsilon$-piece with both $V$ and $\phi(V)$ of length $\frac{1}{9}\|r\|$.
% This implies that $r$ has a $\delta'$-piece with $V$ and $\phi(V)$ of length $>\frac{1}{9}\|r\|-2\epsilon >\frac{1}{10}\|r\|$.
% It follows that $r$ and $\phi(r)$ have a $2\delta'$-piece (and hence an $\epsilon$-piece) of length $>\frac{1}{10}\|r\|$. However, $r,\phi(r)\in \mathcal{R} = \phi(\mathcal{R})$, and we get a contradiction if $\mathcal{R}$ has small enough cancellations.

Similarly one proves the case $\min \{ \|r_{z1} \|, \| r_{z4}\|\} > \lambda\ii(\mu\|r\|+2\delta''+c)$ and the path $p''$ contains $v_6$.

Case 3. $\|r_{z4}\|>\omega - \omega_1 - \omega_2=:\omega_3$ (and similarly for $r_{z1}$, $r_{z2}$ and $r_{z3}$). For small enough $\mu$ and large enough $\rho$ we can assume that $\omega_3/\|r\|$ is arbitrarily close to $1$, and thus we can assume $(\|r\|-\omega_3)+2\delta'' < \lambda \omega_3 - c$.
However, as in Case 2 in the proof of Item \ref{lem: summary of properties of sc quotients - quasiconvex} of Lemma \ref{lem: summary of properties of sc quotients}, there exists a shortcut to $z$ (in $G'$), contradicting the assumption that $z$ is a geodesic.
% This is impossible, as it would contradict the assumption that $z$ was a geodesic in $G/\ngen{\mathcal{R}}$.
\end{proof}

\bibliographystyle{abbrv}
\bibliography{biblio}

% \import{./}{addendum.tex}

\end{document}